%% file: nolip.tex
\documentclass[12pt, USenglish, reqno]{article}

\input{TeX/Preamble.tex}

\newcommand{\TheTitle}{Proximal Gradient Algorithms under Local Lipschitz Gradient Continuity}
\newcommand{\TheSubtitle}{A Convergence and Robustness Analysis of PANOC}

\newcommand{\TheKeywords}{%
	Nonsmooth nonconvex optimization\AND
	locally Lipschitz gradient\AND
	for\-ward-backward splitting\AND
	linesearch methods%
}
\newcommand{\TheAMSsubj}{%
	\amsmscLink{49J52}
	\AND
	\amsmscLink{65K05}
	\AND
	\amsmscLink{90C30}
}
\newcommand{\TheFunding}{%
	A. Themelis acknowledges the support of the Japan Society for the Promotion of Science (JSPS) KAKENHI grant JP21K17710.
}
\newcommand{\TheAffiliationADM}{%
	Universit\"at der Bun\-des\-wehr M\"un\-chen\NL
	Department of Aerospace Engineering\NL
	Institute of Applied Mathematics and Scientific Computing\NL
	Werner-Heisenberg-Weg 39, 85577 Neubiberg, Germany\NL[.]%
	\emailLink{alberto.demarchi@unibw.de}\NL
	\orcidLink{0000-0002-3545-6898}%
}
\newcommand{\TheAffiliationAT}{%
	Kyushu University\NL
	Faculty of Information Science and Electrical Engineering (ISEE)\NL
	744 Motooka, Nishi-ku, 819-0395 Fukuoka, Japan\NL[.]%
	\emailLink{andreas.themelis@ees.kyushu-u.ac.jp}\NL
	\orcidLink{0000-0002-6044-0169}%
}

	\author{%
		Alberto De Marchi\thanks{\TheAffiliationADM}%
		\and
		Andreas Themelis\thanks{\TheAffiliationAT}%
	}%
	\date{}
	\title{\TheTitle\thanks{\TheFunding}\texorpdfstring{\large\\--- }{: }\TheSubtitle\texorpdfstring{ ---}{}}

\begin{document}

		\maketitle
		\begin{abstract}
			\input{TeX/Text/Abstract.tex}
		\end{abstract}
		\keywords{\TheKeywords}
		\subclass{\TheAMSsubj}

	{\def\iter{}%
		\section{Introduction}
			\input{TeX/Text/Introduction.tex}

		\section{Problem Setting and Preliminaries}\label{sec:Setting}
			\input{TeX/Text/Setting.tex}

			\subsection{Notational Conventions}\label{sec:Notation}%
				\input{TeX/Text/Setting/Notation.tex}

			\subsection{Proximal Gradient Iterations}
				\input{TeX/Text/Setting/PG.tex}

			\subsection{Forward-Backward Envelope}\label{sec:FBE}%
				\input{TeX/Text/Setting/FBE.tex}
	}%

		\section{Good and Bad Adaptive Stepsize Selection Rules}\label{sec:Algorithm}
			\input{TeX/Text/Algorithm.tex}

			\subsection{\panoc+: the ``Good'' Adaptive Stepsize Rule}
				\input{TeX/Text/Algorithm/Good.tex}

			\subsection{Failure of ``Bad'' \panoc{} without Globally Lipschitz Gradient}\label{sec:Counterexample}
				\input{TeX/Text/Algorithm/Counterexample.tex}

			\subsection{``Good'' \panoc{+} vs ``Bad'' \panoc{}}\label{sec:CounterexampleBox}
				\input{TeX/Text/Algorithm/CounterexampleBox.tex}

		\section{Algorithmic Analysis under Inexact Proximal Oracles}\label{sec:Analysis}
				\input{TeX/Text/Analysis.tex}

			\subsection{Well Definedness and Convergence Results}\label{sec:Results}%
				\input{TeX/Text/Analysis/Results.tex}

			\subsection{Termination Criteria}\label{sec:Termination}%
				\input{TeX/Text/Analysis/Termination.tex}

			\subsection{Nonmonotone Variant}\label{sec:NM}%
				\input{TeX/Text/Analysis/NM.tex}

			\subsection{Adaptive Proximal Gradient Method}
				\input{TeX/Text/Analysis/PG.tex}

		\section{Conclusions}\label{sec:Conclusions}
			\input{TeX/Text/Conclusions.tex}

		\phantomsection
		\addcontentsline{toc}{section}{References}%
		\bibliographystyle{plain}
		\bibliography{nolip.bib}

\end{document}

%% file: TeX/Preamble.tex
	\PassOptionsToPackage{noend}{algpseudocode}
	\PassOptionsToPackage{nameinlink,capitalize}{cleveref}
	\PassOptionsToPackage{breakspaces}{clevethm}
	\PassOptionsToPackage{inline}{enumitem}
	\PassOptionsToPackage{bookmarksdepth=4}{hyperref}
	\PassOptionsToPackage{dvipsnames}{xcolor}
	\usepackage{tcolorbox}
	\usepackage[
		alglinenumber,
		excludethm=ass,
		thmreset=section,
		eqreset=section,
	]{myPreamble}
	\pgfplotsset{compat=1.16}

	\Alignfalse
	\newcommand{\NL}[1][]{\ifstrempty{#1}{,}{#1}\ }
	\newcommand{\AND}{ \(\cdot\) }
	\newcommand{\keywords}[1]{\par\noindent{\def\and{\unskip,\ }{\bf Keywords. }#1}\par}
	\newcommand{\subclass}[1]{\par\noindent{\def\and{\unskip,\ }{\bf AMS subject classifications. }#1}\par}
	\let\OLDand\and
	\def\and{\texorpdfstring{\OLDand}{, }}%

	\Crefname{figure}{Figure}{Figures}
	
	\theoremstyle{plain}
	\newtheorem*{ass}{Blanket assumption}
	\setlistlabel[proofitemize]{\textbullet}

	\crefname{ALG@line}{step}{steps}

	\newtcolorbox{mybox}[1][]{%
		left=0pt,
		right=0pt,
		top=0pt,
		bottom=0pt,
		colback=MidnightBlue!5,
		colframe=MidnightBlue!10,
		width=\dimexpr\textwidth\relax,
		enlarge left by=0mm,
		boxsep=5pt,
		arc=5pt,outer arc=5pt,
		#1
	}

	\newcommand{\panoc}[1]{%
		\texorpdfstring{%
			\ifstrempty{#1}{%
				\hyperref[alg:PANOC]{PANOC}%
			}{%
				\if+#1%
					\hyperref[alg:PANOC+]{PANOC$^{#1}$}%
				\else
					PANOC$^{#1}$%
				\fi
			}%
		}{PANOC#1}%
	}

	\makeatletter
		\newcommand{\Res}{\@ifnextchar_\Res@sub\Res@no@sub}
		\def\Res@sub_#1{\operatorname{R}_{#1}}
		\def\Res@no@sub{\operatorname{R}_{\FBstepsize}}
		
		\newcommand{\T}{\@ifnextchar_\T@sub\T@no@sub}
		\def\T@sub_#1{\operatorname{T}_{#1}}
		\def\T@no@sub{\operatorname{T}_{\FBstepsize}}

		\newcommand{\iter}{k}

		\renewcommand{\FBstepsize}{\@ifstar\@@FBstepsize\@FBstepsize}
		\newcommand{\@FBstepsize}{\gamma_{\iter}}
		\newcommand{\@@FBstepsize}{\gamma_{\iter-1}}

		\newcommand{\FBE}{\@ifstar\@@FBE\@FBE}
		\newcommand{\@FBE}{\varphi_{\FBstepsize}^{\text{\sc fb}}}
		\newcommand{\@@FBE}{\varphi_{{\FBstepsize*}}^{\text{\sc fb}}}

		\newcommand{\LL}{\@ifstar\@@LL\@LL}
		\newcommand{\@LL}{\@ifnextchar_\@LL@sub\@LL@nosub}
		\newcommand{\@@LL}{\@ifnextchar_\@@LL@sub\@@LL@nosub}
		\newcommand{\@LL@nosub}{\mathcal L_{\nicefrac{1}{\FBstepsize}}}
		\newcommand{\@@LL@nosub}{\mathcal L_{\nicefrac{1}{\FBstepsize*}}}
		\def\@LL@sub_#1{\mathcal L_{#1}}
		\def\@@LL@sub_#1{\mathcal L_{#1}}
	\makeatother

	\newcommand{\emailLink}[1]{\href{mailto:#1}{#1}}
	\newcommand{\orcidLink}[1]{ORCID: \href{https://orcid.org/#1}{#1}}
	\newcommand{\amsmscLink}[1]{\href{http://www.ams.org/mathscinet/msc/msc2020.html?t=#1}{#1}}

	\usepgfplotslibrary{groupplots}%
	\usetikzlibrary{matrix}%
	\definecolor{origcolor}{RGB}{55,126,184}%
	\definecolor{pluscolor}{RGB}{228,26,28}%
	\definecolor{origglipcolor}{RGB}{77,175,74}%
	\definecolor{plusglipcolor}{RGB}{10,10,10}%
	\definecolor{plusglipndircolor}{RGB}{10,10,10}%
	\definecolor{origglipndircolor}{RGB}{10,10,10}%

%% file: TeX/Text/Abstract.tex
Composite optimization offers a powerful modeling tool for a variety of applications and is often numerically solved by means of proximal gradient methods.
In this paper, we consider fully nonconvex composite problems under only local Lipschitz gradient continuity for the smooth part of the objective function.
We investigate an adaptive scheme for PANOC-type methods (Stella et al. in Proceedings of the IEEE 56th CDC, 1939--1944, 2017), namely accelerated linesearch algorithms requiring only the simple oracle of proximal gradient.
While including the classical proximal gradient method, our theoretical results cover a broader class of algorithms and provide convergence guarantees for accelerated methods with possibly inexact computation of the proximal mapping.
These findings have also significant practical impact, as they widen scope and performance of existing, and possibly future, general purpose optimization software that invoke PANOC as inner solver.

%% file: TeX/Text/Introduction.tex
Problems involving the minimization of the sum of a smooth and a nonsmooth function are of interest for a wide variety of applications ranging from optimal and model predictive control (MPC), signal processing, compressed sensing, machine learning, and many others; see, \eg \cite{combettes2011proximal,parikh2014proximal,stathopoulos2016operator} and references therein.
Structured problems can also arise as subproblems within other numerical optimization algorithms, \eg the augmented Lagrangian method (ALM) \cite{rockafellar1976augmented,bertsekas1999constrained,birgin2014alm}.
These use cases often yield nonconvex and large-scale problems and can pose stringent requirements in terms of both computation and memory.

In the last few years, these considerations led to a renewed interest in algorithms of splitting nature \cite{combettes2011proximal,parikh2014proximal} owing to their simple operation oracles and low memory footprint, on top of their amenability to address nonsmooth, possibly nonconvex, constrained  problems, making them widely applicable.
The price of this flexibility is paid in terms of slow convergence and sensitivity to ill conditioning, hindering their direct employment to real-time applications, such as MPC, where optimal solutions to hard problems have to be retrieved in very limited time.
	
Inspired by Newton-type methods for smooth optimization, second-order information can be adopted, so as to better scale with problem size and achieve asymptotic superlinear rates.
However, only local convergence guarantees can be expected without introducing a globalization strategy, such as a backtracking linesearch procedure.
Unfortunately, for nonsmooth problems, even if fast search directions are available classical linesearch strategies are not applicable.
In fact, lacking directional differentiability, the notion of descent directions is not relevant for possibly extended-real-valued, discontinuous functions.

In this very setting, the recently introduced PANOC \cite{stella2017simple} demonstrated how these downsides within the proximal gradient (PG) algorithm can be overcome while retaining all the favorable features.
Essentially, PANOC is a linesearch method that uses the so-called forward-backward envelope (FBE) \cite{patrinos2013proximal} as merit function to globalize the convergence of fast local methods.
It offers an umbrella framework that includes the PG method as special instance; other variations are obtained by selecting virtually arbitrary update directions, which are suitably dampened in such a way to guarantee convergence.
A most prominent use case is the employment of directions stemming from methods of quasi-Newton type applied to the nonlinear equation \(\Res(x)=0\) that encodes first-order necessary conditions for optimality, where \(\Res\) is a (set-valued) generalization of the gradient mapping for nonsmooth problems, cf. \eqref{eq:Res}.
In accommodating arbitrary update directions, PANOC does not require differentiability properties on the merit function and waives the need of regularization terms to enforce a descent condition on the update directions.
We defer a more detailed analysis to the dedicated \cref{sec:Algorithm}.

Although the algorithm uses the same computational oracle of PG, curvature information enables asymptotic superlinear rates under mild assumptions at the limit point \cite{stella2017simple}.
By employing directions of \emph{quasi}-Newton type, no inner iterative procedure nor Hessian evaluations are required, making PANOC's iterations simple, lightweight and scalable.
Because of these favorable properties, PANOC was originally meant as a nonlinear MPC solver particularly suited for embedded applications subject to limited hardware capabilities, such as
land and aerial vehicles \cite{sathya2018embedded,small2019aerial,katriniok2019nonlinear}
and robotics \cite{astudillo2020towards,sathya2020real,berlin2021trajectory};
see also \cite{paalsson2020nonlinear,hermans2021penalty}
for extensive surveys and comparisons with other popular methods.
Its success in the field led to a reconsideration of the spectrum of problems that the solver could be applied to.
On a historical note, this evolution was reflected by a swift rebranding of the acronym over the years, originally meant as \emph{Proximal Averaged Newton-type method for Optimal Control} in the original publication \cite{stella2017simple}, but then tacitly reproposed as the same method \emph{for Optimality Conditions} in \cite{antonello2020proximal} (and subsequent appearances) to allude to its applicability to the much broader range of composite minimization problems.
This flexibility was further exploited in \cite{sopasakis2020open}, where PANOC is employed as inner solver for ALM minimization subproblems for the general purpose Optimization Engine (OpEn) solver.

This rapid evolution was perhaps neglectful of some aspects, primarily because PG is subject to binding assumptions to guarantee a global Lipschitz differentiability requirement.
In the context of MPC, physical bounds on input variables result in optimization problems where the feasible set is bounded, in which case \emph{local} Lipschitzness can be shown to suffice, making virtually no exclusion to the problems that can be addressed.
In more general formulations, and especially so in a fully nonconvex setting, however, all known results are valid under a \emph{global} Lipschitzness assumption, with the very recent work \cite{kanzow2021convergence} possibly emerging as unique exception in a vast literature; see also \cite{cruz2016convergence,salzo2017variable} for convex problems.
Other alternatives are to be found in the Bregman setting \cite{bolte2018first,lu2018relatively,ahookhosh2021bregman}, which are however subject to (and thus limited in applicability by) the identification of a distance-generating function enabling a so-called Lipschitz-like convexity condition and that makes induced proximal operations tractable at the same time.
While this may not seem a major issue in composite minimization, it undeniably constitutes a severe drawback in ALM contexts, where constraints relaxation can produce subproblems with unbounded feasible sets, without this necessarily being the case for the original problem.
Although adding large box constraints to ensure convergence may be thought of as a viable solution, unsatisfactory practical performance can persist because of poor geometry estimation, as we will show.

This paper addresses the above-mentioned shortcomings of PANOC, and of PG as a byproduct, by investigating an adaptive stepsize selection rule for its PG oracle.
This criterion, in a slightly less general form, was first proposed in \cite[Alg. 7]{pas2021matrix}, but without theoretical guarantees and driven from a different observation, namely the poor performance of PANOC if initial stepsizes are badly estimated.
After confirming this claim with case study examples, we provide a complete convergence theory showing that the method, here referred to as \panoc+ for clarity, can also cope with \emph{local} Lipschitzness, while this is not the case for PANOC.
Furthermore, we examine the robustness of the improved method with respect to suboptimal solutions of the PG subproblems.
These findings will significantly impact on \panoc{(+)}, both in performance and applicability, propagating to all its dependencies, \eg by removing stringent assumptions of general purpose optimization solvers such as OpEn \cite{sopasakis2020open}.
Indeed, the significance and effectiveness of \panoc{+} have already been demonstrated in \cite{pas2021alpaqa,demarchi2022constrained}.
As part of the open-source Julia package ProximalAlgorithms.jl \cite{stella2022proximalalgorithms}, our implementation \texttt{PANOCplus} of \panoc{+} is publicly available.

A convergence analysis of PG with a locally Lipschitz smooth term and possibly inexact inner minimizations is obtained as simple byproduct of the more general theory here developed.
Indeed, a vast class of algorithms is covered by the analysis in this work, thanks to the arbitrarity of the selected update directions within the PANOC framework.

%% file: TeX/Text/Setting.tex
In this paper we consider structured minimization problems
\[\tag{P}\label{eq:P}
	\minimize_{x\in\R^n}\varphi(x)\coloneqq f(x)+g(x),
\]
where $x\in\R^n$, $n\in\N$, is the decision variable, under the following standing assumptions, assumed throughout.

\begin{mybox}
	\begin{ass}
		The following hold in problem \eqref{eq:P}:
		\begin{enumeratass}
		\item\label{ass:f}%
			\(\func{f}{\R^n}{\R}\) has locally Lipschitz-continuous gradient.
		\item\label{ass:g}%
			\(\func{g}{\R^n}{\R\cup\set\infty}\) is proper, lsc, and \(\gamma_g\)-prox-bounded.
		\item\label{ass:phi}%
			\(\inf\varphi>-\infty\).
		\end{enumeratass}
	\end{ass}
\end{mybox}

Motivated by its efficiency and popularity, yet aware of its inability to address this general problem formulation, this paper studies a robustified variant of \panoc{} algorithm with adaptive stepsize selection \cite[Rem. III.4]{stella2017simple}, building upon the preliminary work of \cite[\S6.1]{pas2021matrix}.
\panoc{} and the proposed generalization \panoc+ will be presented and compared in \cref{sec:Algorithm}, after the needed definitions and preliminary material are covered in this section.

%% file: TeX/Text/Setting/Notation.tex
With \(\R\) and \(\Rinf\coloneqq\R\cup\set\infty\) we denote the real and extended-real line, and by \(\N=\set{0,1,\dots}\) the set of natural numbers.
The effective domain of an extended-real-valued function \(\func{h}{\R^n}{\Rinf}\) is denoted by \(\dom h\coloneqq\set{x\in\R^n}[h(x)<\infty]\), and we say that \(h\) is: proper if \(\dom h\neq\emptyset\); lower semicontinuous (lsc) if \(h(\bar x)\leq\liminf_{x\to\bar x}h(x)\) for all \(\bar x\in\R^n\); coercive if \(h(x)\to\infty\) as \(\|x\|\to\infty\).
For $\alpha\in\R$, the $\alpha$-\emph{sublevel set} of $h$ is $\lev_{\leq\alpha} h \coloneqq \set{x\in\R^n : h(x) \leq \alpha}$.

The notation \(\ffunc{T}{\R^n}{\R^n}\) indicates a set-valued mapping \(T\) that associates every \(x\in\R^n\) to a subset \(T(x)\subseteq\R^n\).
The \emph{graph} of \(T\) is \(\graph T\coloneqq\set{(x,y)}[y\in T(x)]\).
Following \cite[Def. 8.3]{rockafellar1998variational}, we denote by $\ffunc{\hat{\partial} h}{\R^n}{\R^n}$ the \emph{regular (Fr{\'e}chet) subdifferential} of $h$, where
\begin{equation}
	v \in \hat{\partial} h(\bar{x})
	\quad\defeq[\Leftrightarrow]\quad
	\liminf_{\substack{x\to\bar{x}\\x\neq\bar{x}}}
	\frac{h(x) - h(\bar{x}) - \langle v, x-\bar{x}\rangle}{\|x-\bar{x}\|}
	\geq 0 .
\end{equation}
The \emph{(limiting) subdifferential} of $h$ is $\ffunc{\partial h}{\R^n}{\R^n}$, where $v \in \partial h(\bar{x})$ if there exist sequences $\seq{x^k,v^k}$ in \(\graph\hat\partial f\) such that $(x^k,v^k,h(x^k))\to(\bar x,v,h(\bar x))$.
These subdifferentials of $h$ at $\bar{x}\in\R^n$ satisfy \(\hat\partial(h+h_0)(\bar x)=\hat\partial h(\bar x)+\nabla h_0(\bar x)\) and \(\partial(h+h_0)(\bar x)=\partial h(\bar x)+\nabla h_0(\bar x)\) for any $\func{h_0}{\R^n}{\Rinf}$ continuously differentiable around $\bar{x}$ \cite[Ex. 8.8]{rockafellar1998variational}.
With respect to \eqref{eq:P}, we say that \(x^*\in\dom\varphi\) is \emph{stationary} if \(0\in\partial\varphi(x^*)\), which constitutes a necessary optimality condition of \(x^*\) for the minimization of \(\varphi\) \cite[Thm. 10.1]{rockafellar1998variational}.

Given a parameter value \(\gamma>0\), the \emph{Moreau envelope} function \(h^\gamma\) and the \emph{proximal mapping} \(\prox_{\gamma h}\) are defined by
\begin{align}
	\label{eq:moreauenvelope}
	h^\gamma(x)
	{}\coloneqq{} &
	\inf_{z\in\R^n}\set{
		h(z)
		{}+{}
		\tfrac{1}{2\gamma}\|z-x\|^2
	} ,
	\\
	\label{eq:proxoperator}
	\prox_{\gamma h}(x)
	{}\coloneqq{} &
	\argmin_{z\in\R^n}\set{
		h(z)
		{}+{}
		\tfrac{1}{2\gamma}\|z-x\|^2
	},
\end{align}
and we say that \(h\) is \emph{prox-bounded} if it is proper and \(h + \tfrac{1}{2\gamma}\|\cdot\|^2\) is bounded below on \(\R^n\) for some \(\gamma>0\).
The supremum of all such \(\gamma\) is the threshold \(\gamma_h\) of prox-boundedness for \(h\).
In particular, if \(h\) is bounded below by an affine function, then \(\gamma_h = \infty\).
When \(h\) is lsc, for any \(\gamma\in(0,\gamma_h)\) the proximal mapping \(\prox_{\gamma h}\) is nonempty- and compact-valued, and the Moreau envelope \(h^\gamma\) finite and locally Lipschitz continuous \cite[Thm. 1.25 and Ex. 10.32]{rockafellar1998variational}.

%% file: TeX/Text/Setting/PG.tex
Given a point \(x \in \R^n\), one iteration of the proximal gradient (PG) method for problem \eqref{eq:P} consists in selecting
\begin{equation}\label{eq:FBSiter}
	\def\iter{}
	\bar x
	{}\in{}
	\T(x)
	{}\coloneqq{}
	\FB{x},
\end{equation}
where \(\gamma\in(0,\gamma_g)\) is a stepsize parameter.
The necessary optimality condition in the minimization problem defining the proximal mapping then reads
\begin{equation}\label{eq:proxoptim}
	\tfrac{1}{\gamma}(x-\bar x)
	{}-{}
	(\nabla f(x)-\nabla f(\bar x))
{}\in{}
	\hat\partial\varphi(\bar x),
\end{equation}
and in particular the fixed-point inclusion \(x\in\T(x)\) implies the stationarity condition \(0\in\partial\varphi(x)\).
By interpreting \eqref{eq:FBSiter} as a fixed-point iteration, one can also consider the associated (set-valued) fixed-point residual \(\Res\), namely
\begin{equation}\label{eq:Res}
	\Res(x)
	{}\coloneqq{}
	\tfrac{1}{\gamma}\bigl(x-\T(x)\bigr),
\end{equation}
and seek fixed points of \(\T\) as zeros of the residual \(\Res_\gamma\).

%% file: TeX/Text/Setting/FBE.tex
At the heart of PANOC rationale is the observation that, under assumptions, the fixed-point residual \(\Res\) in \eqref{eq:Res} is continuous around and even differentiable at critical points \cite[\S 4]{themelis2018forward}, and the inclusion problem \(0\in\Res({}\cdot{})\) reduces to a well-behaved system of equations, when close to solutions.
This motivated the adoption of Newton-type directions on \(\Res\), that enable fast convergence when close to solutions.
The key tool enabling convergence regardless of whether or not the initial point happens to be sufficiently close to a solution is the so-called forward-backward envelope (FBE).

\begin{defin}[forward-backward envelope]%
	\begin{subequations}\label{subeq:FBE}%
		Relative to \eqref{eq:P}, the FBE with stepsize \(\gamma\in(0,\gamma_g)\) is
		\begin{align}
		\label{eq:FBE}
			\FBE(x)
		{}\coloneqq{} &
			\min_{w\in\R^n}\set{
				f(x)
				{}+{}
				\innprod{\nabla f(x)}{w-x}
				{}+{}
				g(w)
				{}+{}
				\tfrac{1}{2\gamma}\|w-x\|^2
			}
		\\
		\label{eq:FBEMoreau}
		{}={} &
			f(x)
			{}-{}
			\tfrac\gamma2\|\nabla f(x)\|^2
			{}+{}
			g^\gamma(\Fw x)
		\shortintertext{%
			or, equivalently, letting \(\bar x\) be any element of \(\T(x)\),
		}
		\label{eq:FBEmin}
		{}={} &
			f(x)
			{}+{}
			\innprod{\nabla f(x)}{\bar x-x}
			{}+{}
			g(\bar x)
			{}+{}
			\tfrac{1}{2\gamma}\|\bar x-x\|^2.
		\end{align}
	\end{subequations}
\end{defin}

Owing to its continuity properties, the FBE has been employed to generalize and improve PG-based algorithms that address the general setting of structured nonconvex optimization \cite{liu2017further,themelis2018forward,bonettini2020convergence}.
The following results are well known when \(f\) has globally Lipschitz gradient \cite[Prop.s 4.2 and 4.3]{themelis2018forward}.
A simple proof in the more general setting addressed here is given for completeness.

\begin{lem}[Properties of the FBE]\label{thm:FBE}%
	\def\iter{}%
	For any \(\gamma\in(0,\gamma_g)\) the following hold:%
	\begin{enumerate}
	\item\label{thm:FBE:C0}%
		\(\FBE\) is real valued and strictly continuous.
	\item\label{thm:FBE:leq}%
		\(\FBE(x)\leq\varphi(x)\) for any \(x\in\R^n\), with equality holding iff \(x\in\T(x)\).
	\item\label{thm:FBE:geq}%
		If \(\bar x\in\T(x)\) and
		\(
			f(\bar x)
		{}\leq{}
			f(x)
			{}+{}
			\innprod{\nabla f(x)}{\bar x-x}
			{}+{}
			\tfrac L2\|\bar x-x\|^2
		\),
		then
		\begin{equation}\label{eq:FBE:geq}
			\FBE(\bar x)
		{}\leq{}
			\varphi(\bar x)
		{}\leq{}
			\FBE(x)
			{}-{}
			\tfrac{1-\gamma L}{2\gamma}\|x-\bar x\|^2.
		\end{equation}
	\end{enumerate}
	\begin{proof}
		Assertion \ref{thm:FBE:C0} follows from the expression \eqref{eq:FBEMoreau}, owing to the similar property of the Moreau envelope \(g^\gamma\), while \ref{thm:FBE:leq} is obtained by taking \(w=x\) in \eqref{eq:FBE}.
		The first inequality in \ref{thm:FBE:geq} owes to item \ref{thm:FBE:leq} (independently of \(L\)), and the second one follows from the expression \eqref{eq:FBEmin} of \(\FBE\).
	\end{proof}
\end{lem}

%% file: TeX/Text/Algorithm.tex
{\def\iter{}%
	As briefly mentioned in \cref{sec:FBE}, the FBE is the key tool for \emph{globalizing} the convergence of fast local methods, such as of quasi-Newton type, applied to the nonlinear equation \(\Res(x)=0\) encoding necessary optimality conditions for \eqref{eq:P}.
	Elaborating on how Newton-type directions can be selected given the nonsmooth, possibly set-valued, nature of \(\Res\) is beyond the scope of this survey, and the interested reader is referred to \cite{themelis2018forward,stella2017simple}.
	The core idea is nevertheless the same as in the familiar context of smooth minimization: trying to enforce (supposedly fast) updates \(x\mapsto x+d\) in place of ``nominal'' updates \(x\mapsto\bar x\), where \(\bar x\) would amount to a gradient step or, in our nonsmooth setting, a proximal gradient step \(\bar x\in\T(x)\) as in \eqref{eq:FBSiter}.
	Still in complete analogy with the smooth case, accepting a candidate update \(x+d\) must be validated by a ``quality check'', like an Armijo-type condition, in violation of which \(d\) is either discarded or dampened with a smaller stepsize.
	\panoc{} is precisely a mechanism to dampen and accept update directions in a nonsmooth setting, using the FBE as validation control.
	Its steps are given in \cref{alg:PANOC}.

	\begin{algorithm}[t]%
		\caption{Original PANOC with ``bad'' adaptive stepsize \(\gamma\) \cite[Rem. III.4]{stella2017simple}}%
		\label{alg:PANOC}%
		\input{TeX/Alg/PANOC.tex}	\end{algorithm}

	A basic assumption for \panoc{} is that \(\nabla f\) be globally \(L_f\)-Lipschitz, so that a well-known quadratic upper bound, see \eg \cite[Prop. A.24]{bertsekas1999nonlinear}, ensures that \(L=L_f\) can be taken for all \(x\in\R^n\) in \cref{thm:FBE:geq}.
	Alternatively, if $g$ has bounded domain and the selected directions $d^k$ are bounded, it suffices that $\nabla f$ is locally Lipschitz-continuous; see \cite[Rem. III.4]{stella2017simple}.
	For any \(\alpha\in(0,1)\) the choice \(\gamma_k=\nicefrac{\alpha}{L_f}\) then violates \cref{state:PANOC:gammaLS}, meaning that \(\gamma_k\equiv\gamma\) is constant.
	The dampening of the direction occurs at \cref{state:PANOC:x+}, where starting with \(\tau_k=1\) the candidate update \(x^{k-1}+d^k\) is pushed towards \(\bar x^{k-1}\in\T(x^{k-1})\) by reducing the steplength \(\tau_k\) until the value of the FBE is sufficiently reduced, cf. \cref{state:PANOC:tauLS}.
	The process terminates, since \(\FBE\) is continuous (at \(\bar x^{k-1}\)), and it is strictly smaller than
	\(
		\FBE(x^{k-1})
		{}-{}
		\beta\tfrac{1-\alpha}{2\gamma_{k-1}}
		\|\bar x^{k-1}-x^{k-1}\|^2
	\)
	there, cf. \eqref{eq:FBE:geq}.

%% file: TeX/Alg/PANOC.tex
\begin{algorithmic}[1]%
	\def\iter{k-1}%
	\linespread{1.4}\selectfont%
	\Require
		\(x^0\in\R^n\);~~
		\(\gamma_0\in\bigl(0,\gamma_g\bigr)\);~~
		\(D\geq0\);~~
		\(\alpha,\beta\in(0,1)\)%
	\Initialize
		\def\iter{0}%
		\(k=0\),~
		compute \(\bar x^0\in\T(x^0)\),~
		and start from \cref{state:PANOC:first}
\vspace*{1ex}%
\hrule
\vspace*{1ex}%
\def\iter{k-1}%
	\State
		\label{state:PANOC:init}%
		Select an update direction \(d^k\in\R^n\) with \(\|d^k\|\leq D\|\bar x^{k-1}-x^{k-1}\|\) and set \(\tau_k=1\)%
	\State
		\label{state:PANOC:x+}%
		\(x^k=(1-\tau_k)\bar x^{k-1}+\tau_k(x^{k-1}+d^k)\)
	\State
		\label{state:PANOC:barx}%
		Compute \(\bar x^k\in\T(x^k)\) and use it to evaluate \(\FBE(x^k)\) as in \eqref{eq:FBEmin}
	\If{ \(
		\FBE(x^k)
	{}>{}
		\FBE(x^{k-1})
		{}-{}
		\beta\tfrac{1-\alpha}{2\gamma_{k-1}}
		\|\bar x^{k-1}-x^{k-1}\|^2
	\) }%
	\label{state:PANOC:tauLS}%
		\Statex*
			\(\tau_k\gets\nicefrac{\tau_k}{2}\)~
			and go back to \cref{state:PANOC:x+}
	\EndIf
\def\iter{k}%
	\State
		\(\gamma_k\gets\gamma_{k-1}\)
	\While{ \(f(\bar x^k)>f(x^k)+\innprod{\nabla f(x^k)}{\bar x^k-x^k}+\tfrac{\alpha}{2\gamma_k}\|\bar x^k-x^k\|^2\) }
		\label{state:PANOC:first}\label{state:PANOC:gammaLS}%
		\Statex*
			\(\gamma_k\gets\nicefrac{\gamma_k}{2}\)
			~and recompute~
			\(\bar x^k\in\T(x^k)\)
	\EndWhile
	\State
		\(k\gets k+1\)~
		and start the next iteration at \cref{state:PANOC:init}
\end{algorithmic}

%% file: TeX/Text/Algorithm/Good.tex
What is presented in \cref{alg:PANOC} is actually the ``adaptive'' variant of \panoc{}, which still works under the assumption of global Lipschitz differentiability but waives the need of prior knowledge about \(L_f\).
The \(\gamma\)-backtracking at \cref{state:PANOC:gammaLS} decreases (\ie, ``adapts'') \(\gamma_k\) and terminates as soon as the needed bound as in \cref{thm:FBE:geq} is satisfied.
As first noted in \cite[\S6.1]{pas2021matrix}, however, this adaptive criterion may produce bad estimates of the local Lipschitz constant of \(\nabla f\) and overall result in poor algorithmic performance.
The phenomenon can be attributed to an asynchrony between the two backtracking steps, the one dampening the update direction and the one adaptively adjusting the proximal gradient stepsize.
This claim can be verified in the iteration mismatch between variable \(x^k\) and stepsize \(\gamma_{k-1}\) occurring at \cref{state:PANOC:barx} (cf. \cref{rem:notation}).

To account for this fact, \cite[Alg. 7]{pas2021matrix} proposes to adapt the PG stepsize \(\gamma_k\) within the linesearch on the update direction.
As recently showcased in \cite{pas2021alpaqa}, not only does this conservatism prove beneficial in preventing the acceptance of poor quality directions, but it often also reduces the overall computational cost.
Although numerical simulations indicate superior performance, this refined linesearch lacks a theoretical analysis of its convergence properties.

\begin{algorithm}[tb!]
	\caption{\panoc+: the ``good'' adaptive \(\gamma\)-stepsize rule}%
	\label{alg:PANOC+}%
	\input{TeX/Alg/PANOC+.tex}\end{algorithm}

This modification, which we allusively call the ``good'' adaptive variant (or \panoc+ for brevity), is depicted in \cref{alg:PANOC+}.
In fact, the method presented here presents a slight, but important generalization, namely in allowing the selection of a new direction \(d^k\) every time the stepsize \(\gamma_k\) is reduced, cf. \cref{state:PANOC+:gammaLS}, which was not considered in \cite[Alg. 7]{pas2021matrix}.
This flexibility is crucial: whenever the stepsize \(\gamma_k\) changes so does the PG residual mapping \(\Res\), and consistently so should directions using its curvature information.
Moreover, we provide theoretical guarantees on the finite termination of the backtracking linesearch procedure, even without global Lipschitz gradient continuity and merely suboptimal proximal computation.
These findings uphold the algorithmic framework proposed in \cite{stella2017simple,pas2021matrix,pas2021alpaqa} on two aspects: the adaptive linesearch is shown to terminate, and can cope with a merely locally Lipschitz-differentiable term \(f\).
These findings are of high significance also for other methods that rely on PANOC as internal solver, such as the general purpose OpEn \cite{sopasakis2020open}.
Moreover, it will be shown that all this remains true even if the minimization problem defining the PG mapping \(\T\) is solved inexactly and/or suboptimally.

The peculiarity of \panoc+ over the \emph{bad} adaptive rule of original \panoc{} is that the two backtracking steps, the one on the direction \(\tau_k\) and the one on the PG stepsize \(\gamma_k\), are tightly intertwined.
The intricate structure emerges at \cref{state:PANOC+:gammaLS,state:PANOC+:tauLS}: the direction stepsize \(\tau_k\) resets every time the proximal stepsize \(\gamma_k\) is adjusted and, conversely, the value of \(\gamma_k\) is assessed anew when \(\tau_k\) changes.
This entanglement allows the evaluation of the FBE at \cref{state:PANOC+:barx} with an up-to-date stepsize \(\gamma_k\), as opposed to (and eliminating) the asynchrony obstructing \panoc{}'s performance.
The adaptivity of \panoc+ allows the FBE \(\FBE\) to better capture the (local) landscape of \(\varphi\) and, ultimately, to relax the assumption of globally Lipschitz gradient.

To substantiate these claims, in the following \cref{sec:Counterexample} we first showcase the ineffectiveness of \panoc{} applied to problem \eqref{eq:P} where \(f\) has only locally Lipschitz-con\-tin\-u\-ous gradient, and then compare the ``good'' and the ``bad'' adaptive strategies on a common ground in \cref{sec:CounterexampleBox}.

\begin{rem}[Algorithm notation]\label{rem:notation}%
	\Cref{alg:PANOC+} operates two linesearch steps within each iteration, one on the ``proximal'' stepsize \(\gamma_k\) at \cref{state:PANOC+:gammaLS} and one on the ``direction'' stepsize \(\tau_k\) at \cref{state:PANOC+:tauLS}.
	Whenever the respective needed conditions are violated, either \(\gamma_k\) or \(\tau_k\) is reduced and the iteration restarted from a previous step.
	As a consequence, variables may be \emph{overwritten} within each iteration before being accepted.
	To avoid a heavy double-index notation, used only within proofs out of full rigor, the sub- and superscript notation is designed to differentiate temporary and permanent variables; specifically, within iteration \(k\) only variables indexed with \(k\) are updated, whereas those indexed with \(k-1\) remain untouched.
	Similar considerations apply to \cref{alg:PANOC}.
\end{rem}

%% file: TeX/Alg/PANOC+.tex
\begin{algorithmic}[1]%
\linespread{1.4}\selectfont%
\Require
	\(x^0\in\R^n\);~~
	\(\gamma_0\in(0,\gamma_g)\);~~
	\(D\geq0\);~~
	\(\alpha,\beta\in(0,1)\)%
\Initialize
	\(k\gets0\),~~
	and start from \cref{state:PANOC+:barx}
\vspace*{1ex}%
\hrule
\vspace*{1ex}%
\def\iter{k}%
	\State
		\label{state:PANOC+:init}%
		\(\gamma_k\gets\gamma_{k-1}\)
	\State
		\label{state:PANOC+:d}%
		Select an update direction \(d^k\in\R^n\) with \(\|d^k\|\leq D\|\bar x^{k-1}-x^{k-1}\|\) and set \(\tau_k=1\)%
	\State
		\label{state:PANOC+:x+}%
		\(x^k=(1-\tau_k)\bar x^{k-1}+\tau_k(x^{k-1}+d^k)\)
	\State
		\label{state:PANOC+:barx}%
		Compute \(\bar x^k\in\T(x^k)\)
		and use it to evaluate \(\Phi_k\coloneqq\FBE(x^k)\) as in \eqref{eq:FBEmin}%
	\If{
		\(
			f(\bar x^k)
		{}>{}
			f(x^k)
			{}+{}
			\innprod{\nabla f(x^k)}{\bar x^k-x^k}
			{}+{}
			\tfrac{\alpha}{2\gamma_k}\|\bar x^k-x^k\|^2
		\)
	}%
	\label{state:PANOC+:gammaLS}%
		\Statex*%
			\(\gamma_k\gets\nicefrac{\gamma_k}{2}\),~
			and go back to \cref{state:PANOC+:d} if \(k>0\), or \cref{state:PANOC+:barx} if \(k=0\)
	\EndIf
	\If{ \(k>0\) ~{\sc and}~ \(
		\Phi_k
		{}>{}
		\Phi_{k-1}
		{}-{}
		\beta\tfrac{1-\alpha}{2\gamma_{k-1}}
		\|\bar x^{k-1}-x^{k-1}\|^2
	\) }\label{state:PANOC+:tauLS}%
		\Statex*
			\(\tau_k\gets\nicefrac{\tau_k}{2}\)~
			and go back to \cref{state:PANOC+:x+}
	\EndIf
	\State
		\(k\gets k+1\)~
		and start the next iteration at \cref{state:PANOC+:init}
\end{algorithmic}

%% file: TeX/Text/Algorithm/Counterexample.tex
Let us consider the minimization of the convex, twice continuously differentiable, coercive function \(\varphi=f+g\), where \(f(x)=\tfrac29|x|^3\) and \(g=0\), namely
\begin{equation}\label{eq:PANOC_Ex_P}
	\minimize_{x\in\R}{
		\varphi(x)
	{}\coloneqq{}
		\tfrac29|x|^3
		{}+{}
		0
	},
\end{equation}
and adopt \panoc{} as given in \cref{alg:PANOC}.
In particular, we choose the directions as
\begin{equation}\label{eq:dk}
	d_k=\tfrac{9}{2\gamma_{k-1}x_{k-1}}(x_{k-1}-\bar x_{k-1}).
\end{equation}
As we are about to show, starting from any \(x_0>0\) this particular choice of directions complies with the bound \(\|d_k\|\leq D\|x_{k-1}-\bar x_{k-1}\|\) for \(D=18\) and satisfies the \(\tau\)-linesearch with \(\tau_k=1\) for every \(k\).
Moreover, the choice \(\alpha=\nicefrac{16}{27}\) leads to a conveniently simple expression for the \(\gamma\)-linesearch, namely
\(
	\gamma_k
{}\leq{}
	\tfrac{1}{2x_k}
\).
As a result, starting from \(x_0>0\) with \(\gamma_0>\tfrac{1}{4x_0}\), the algorithm reduces iterating the following lines
\begin{equation}\label{eq:PANOC_Ex}
	\begin{cases}
		\text{halven \(\gamma_k\) until \(\gamma_k\leq\frac{1}{2x_k}\)}
	\\[3pt]
		\fillwidthof[c]{x_{k+1}}{\bar x_k}
	{}={}
		x_k(1-\frac23\gamma_kx_k)
	\\[3pt]
		x_{k+1}
	{}={}
		x_k
		{}+{}
		\tfrac{9}{2\gamma_kx_k}(x_k-\bar x_k)
	{}={}
		4x_k
	\end{cases}
\end{equation}
and thus produces a sequence \(x_k=x_04^k\) that is diverging, and causes the cost to increase unboundedly.
We now show the claims one by one.
To this end, denoting \(y_k\coloneqq\gamma_kx_k\) throughout, observe that
\begin{equation}\label{eq:PANOC_Ex:FBE}
	\bar x_k
{}={}
	x_k\left(1-\tfrac23|y_k|\right)
\quad\text{and}\quad
	\FBE(x)
{}={}
	\tfrac29|x|^3(1-\gamma_kx).
\end{equation}
\begin{proofitemize}
\item{\it Linesearch on \(\gamma\).}
	For \(x_k>0\) the backtracking on \(\gamma_k\) at \cref{state:PANOC:gammaLS} (after removing a \(\frac29x_k^3\) factor) terminates when
	\begin{equation}\label{eq:PANOC_Ex:gammaLS}
		\left|1-\tfrac23y_k\right|^3
	{}\leq{}
		1-2y_k+\alpha y_k.
	\end{equation}
	To simplify the computation, observe that necessarily \(y_k\leq1\) for inequality \eqref{eq:PANOC_Ex:gammaLS} to hold, and in particular the argument of the absolute value is necessarily positive: in fact, since \(y_k=\gamma_kx_k>0\) and \(\alpha<1\), \eqref{eq:PANOC_Ex:gammaLS} implies
	\(
		\left|1-\tfrac23y_k\right|^3
	{}\leq{}
		1-y_k
	\),
	hence \(y_k\leq 1\).
	After this simplification and by restricting the analysis to \(y_k=\gamma_kx_k>0\), it can be seen that \eqref{eq:PANOC_Ex:gammaLS} has solution
	\(
		0
	{}<{}
		\gamma_k
	{}\leq{}
		\frac{9}{4x_k}\left(
			1
			{}-{}
			\sqrt{1-\tfrac23\alpha}
		\right)
	\).
	For \(\alpha=\nicefrac{16}{27}\), this bound simplifies to
	\(
		0
	{}<{}
		\gamma_k
	{}\leq{}
		\tfrac{1}{2x_k}
	\)
	as claimed.
	This shows the validity of the first line in \eqref{eq:PANOC_Ex}.
	Since \(\gamma_k\) is halvened (only) until it enters this range, one also has that
	\begin{equation}\label{eq:PANOC_Ex:gamma}
		y_k\coloneqq\gamma_kx_k>\tfrac14
	\quad
		\forall k.
	\end{equation}
\item
	{\it Bound on the directions \(\|d_{k+1}\|\leq D\|x_k-\bar x_k\|\)}.
	Since
	\(
		d_{k+1}=\frac{9}{2\gamma_kx_k}(x_k-\bar x_k)
	\),
	one has
	\(
		\|d_{k+1}\|
	{}={}
		\frac{9}{2|\gamma_kx_k|}\|x_k-\bar x_k\|
	{}\leq{}
		18\|x_k-\bar x_k\|
	\)
	as it follows from \eqref{eq:PANOC_Ex:gamma}.
\item{\it Linesearch on \(\tau\).}
	Starting with \(x_k>0\) we show that
	\(
		x_{k+1}
	{}={}
		x_k+d_{k+1}
	{}={}
		4x_k
	\)
	satisfies the linesearch condition.
	Indeed, by using the expression for the FBE in \eqref{eq:PANOC_Ex:FBE}, according to \cref{state:PANOC:tauLS} the iterate \(x_{k+1}=4x_k\) is accepted if
	\[
		\tfrac29(4x_k)^3
		(1-4y_k)
	{}\leq{}
		\tfrac29x_k^3(1-y_k)
		{}-{}
		\beta(1-\alpha)\tfrac29x_k^3y_k
	\]
	which is easily reduced to
	\(
		y_k
	{}\geq{}
		\frac{
			4^3
			{}-{}
			1
		}{
			4^4
			{}-{}
			1
			{}-{}
			\beta(1-\alpha)
		}
	\).
	Since \(\beta(1-\alpha)<1\), one has
	\(
		\frac{
			4^3
			{}-{}
			1
		}{
			4^4
			{}-{}
			1
			{}-{}
			\beta(1-\alpha)
		}
	{}\leq{}
		\frac{
			4^3
			{}-{}
			1
		}{
			4^4
			{}-{}
			2
		}
	{}<{}
		\tfrac14
	\),
	and \eqref{eq:PANOC_Ex:gamma} implies that the inequality always holds.
\end{proofitemize}

We stressed that, although we consider an exemplary problem designed to yield simple computations, similar arguments would still apply for \(C^\infty\), strongly convex formulations, \eg \(x^4+x^2\); see also \cref{rem:0}.

%% file: TeX/Text/Algorithm/CounterexampleBox.tex
\subsubsection{Robustness against poor directions}
	In spite of the breakdown demonstrated in \cref{sec:Counterexample}, global convergence guarantees for \panoc{} can be recovered by adding a term \(g\) with bounded domain, as is the case of a possibly large but bounded box constraint, and selecting update directions \(d_k\) that are bounded, see \cite[Rem. III.4]{stella2017simple}.
	Nonetheless, as noted in \cite[\S6.1]{pas2021matrix}, this would scarcely help in practice: early iterations would be agnostic to the large box and exhibit the same diverging behavior until the boundary is approached, at which point a drastically reduced stepsize \(\gamma\) would be the cause of a painfully slow convergence.

\begin{figure}[b!]
		\centering
		\includetikz[width=.9\linewidth]{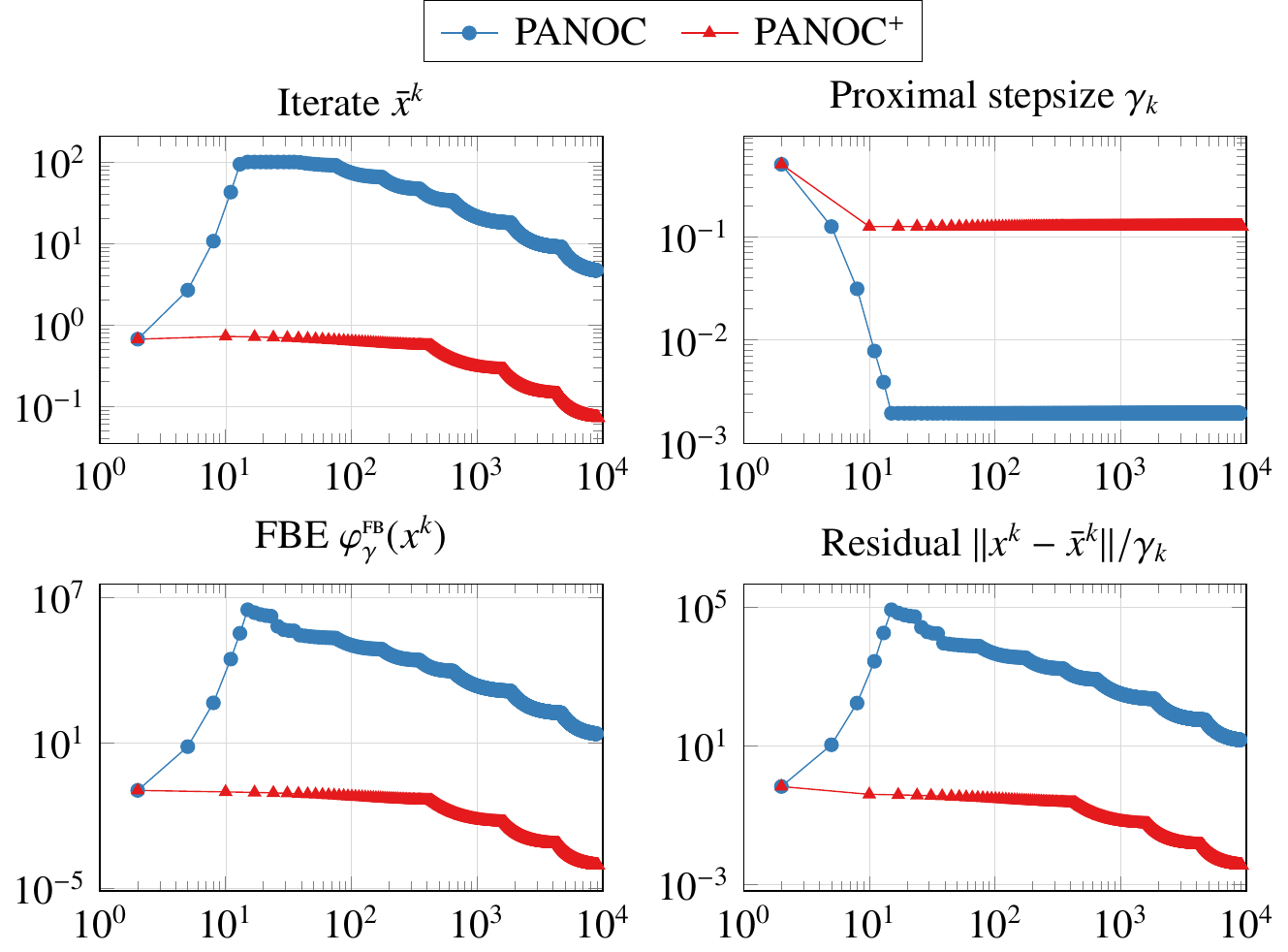}%
		\caption{%
			Comparison of convergence metrics vs number of evaluations of \(\operatorname{T}_\gamma\) for \panoc{} and \panoc+ on the illustrative problem \eqref{eq:PANOC_Ex_P_box}, with directions as in \eqref{eq:dk} saturated in the interval \([-100,100]\).
			We used \(x^0 = 1\), \(\gamma_0 = 1\), \(\alpha = 0.95\), and \(\beta = 0.5\).
			\panoc{}'s iterates diverge until the (safeguarding) box constraint activates, and only then, with a reduced stepsize \(\gamma\), slowly recovers.%
		}%
		\label{fig:counterexamplebox}%
	\end{figure}

	We substantiate these claims by considering the example in \cref{sec:Counterexample} with some amendments.
	In particular, we let \(g\) be the indicator function of the interval \([-B,B]\), namely \(g(x)=0\) if \(|x|\leq B\) and \(g(x)=\infty\) otherwise, and select directions \(d_k\) as above if \(\|d_k\|\leq E\) and \(Ed_k/\|d_k\|\) otherwise, with possibly large but bounded \(B,E\geq 0\).
	The problem becomes
	\begin{equation}\label{eq:PANOC_Ex_P_box}
		\minimize_{x\in\R}\tfrac29|x|^3
	\quad\stt{}
		|x|\leq B.
	\end{equation}
	Adopting these precautions, \panoc{} generates iterates that converge to a solution, starting from any initial point.
	We set \(B=E=100\) for the results displayed in \cref{fig:counterexamplebox} with a comparison against \panoc+.
	Although the latter solves the illustrative problem in its original form (that is, with \(B=\infty\)), we stress that it would not be affected by the safeguards put in place to guarantee the convergence of ``bad'' \panoc{}.

\begin{figure}[t!]
		\centering
		\includetikz[width=.5\linewidth]{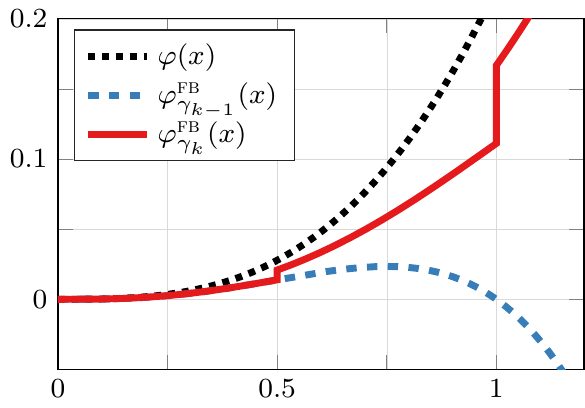}%
		\caption{%
			Comparison of the cost function \(\varphi\) for the illustrative problem \eqref{eq:PANOC_Ex_P} against \panoc{}'s and \panoc+'s merit functions with previous, or initial, estimate \(\gamma_{k-1}=1\)%
		}%
		\label{fig:counterexamplebox_merit}%
	\end{figure}

	The diverging behavior of \panoc{} is apparent, until the safeguards activate, as expected from \cref{sec:Counterexample}.
	At \cref{state:PANOC:barx} \panoc{} accepts an update \(x^k\) based on the sufficient decrease of a merit function defined by the FBE with the \emph{previous} stepsize \(\gamma_{k-1}\).
	\Cref{fig:counterexamplebox_merit} illustrates this phenomenon by comparing the merit functions adopted by \panoc{} and \panoc+ to verify whether a tentative update is to be accepted or not.
	In this example, \panoc{}'s merit function are lower unbounded (see \eqref{eq:PANOC_Ex:FBE}) and full steps along the update directions \(d_k\) are accepted, in fact \emph{favored}, leading to diverging iterates.
	In turn, this results in a temporary departure from the solution, degrading the overall efficiency of the algorithm.
	Conversely, at \cref{state:PANOC+:barx} \panoc+ verifies sufficient decrease of the FBE with the \emph{current} stepsize \(\gamma_k\), yielding monotone decrease of the (time varying, but lower bounded) merit function \(\FBE\), as depicted in \cref{fig:counterexamplebox}.
	Note that the merit function for \panoc+ in \cref{fig:counterexamplebox_merit} is only piecewise continuous because its evaluation is always preceded by the \(\gamma\)-stepsize backtracking, \ie, the stepsize \(\gamma_k=\gamma_k(x^k)\) in \(\FBE\) depends on the candidate update \(x^k\) being tested.
	This adaptivity allows \panoc+ to well estimate the geometry of the cost function \(\varphi\) and to construct a tighter merit function.

	These simulations also show that, despite the more conservative linesearch, \panoc+ does not necessarily require more iterations nor function evaluations to provide a more consistent performance, nor does it lead to a smaller stepsize.
	Indeed, considering larger box constraints and update directions, \ie larger values for \(B\), the limitations and inadequacy of ``bad'' \panoc{} in this setting become apparent, while providing support in favor of the (initially) more conservative adaptive scheme of ``good'' \panoc+.

	\begin{rem}\label{rem:0}%
		Noticeably, the ``bad'' \panoc{} can exhibit this diverging behavior even when the problem admits just one feasible point.
		To see this, let us consider once again the illustrative example above with \(B=0\), so that \(\dom g=\dom\varphi=\set0\).
		Then, patterning the proof in \cref{sec:Counterexample}, we obtain that the algorithm produces a sequence \(\seq{x_k}\) that is diverging, despite the fact that \(\bar x^k=0\) for every \(k\),
		since
		\(
			\FBE(x)
			{}={}
			x^2\bigl(\tfrac{1}{2\gamma_k} - \tfrac49|x|\bigr)
		\)
		is still lower unbounded for any \(\gamma_k>0\).
		This also confirms the necessity of imposing bounded \(\|d^k\|\) in \cite[Rem. III.4]{stella2017simple},
		in addition to \(\|d^k\|\leq D\|x^{k-1}-\bar x^{k-1}\|\) as in \cref{state:PANOC:init}, not needed in the ``good'' \panoc+ even with unbounded domains.
	\end{rem}

\subsubsection{Robustness against poor initial stepsize estimation}
	The poor performance of \panoc{} on problem \eqref{eq:PANOC_Ex_P_box} can be attributed to the bad quality of update directions \(d^k\).
	We now consider a more meaningful comparison on problem \eqref{eq:PANOC_Ex_P_box}, this time with directions given by a classical Newton-type approach.
	We extend \(f\) linearly outside of the box \([-B,B]\) so as to make it (convex and) globally Lipschitz differentiable without affecting the problem.
	We thus consider
	\begin{equation}\label{eq:PANOC_Ex_P_box_smooth}
		\minimize_{x\in\R}f(x)
	\quad\stt{}
		|x|\leq B,
	\end{equation}
	where
	\[
		f(x)
	{}={}
		\begin{ifcases}
			\tfrac29|x|^3 & |x|\leq B \\[3pt]
			\tfrac23B^2(|x|-\tfrac23B)\otherwise.
		\end{ifcases}
	\]
	Because of the constraints, the problem is nonsmooth.
	Nevertheless, since \(f\) is globally \(L_f\)-Lipschitz differentiable (with \(L_f=\frac23B^2\)), the minimization of \(f\) is equivalent to that of \(\FBE\), when \(\gamma<\nicefrac{1}{L_f}\).
	As such, in the spirit of \cite{themelis2019acceleration} we may select update directions based on a Newton method on the FBE.
	We simulate the scenario in which \(L_f\) is unknown, thereby selecting an initial stepsize \(\gamma_0\) larger than \(\nicefrac{1}{L_f}\).
	Since the cost function is coercive and has a unique stationary point, both methods are guaranteed to converge to the unique solution \(x^\star=0\).

	We consider classical Newton directions
	\begin{equation}\label{eq:PANOC_Ex_P_box_smooth_FBEdirection}
		\def\iter{k}d^k=-\max\set{\mu,\,\nabla^2\FBE(x^k)}^{-1}\nabla\FBE(x^k)
	\end{equation}
	with \(\mu>0\) as regularization parameter.
	When not defined, \(\def\iter{k}\nabla^2\FBE\) is intended in a Clarke generalized sense.
	
	\begin{figure}[tb!]%
		\includetikz[width=0.9\linewidth]{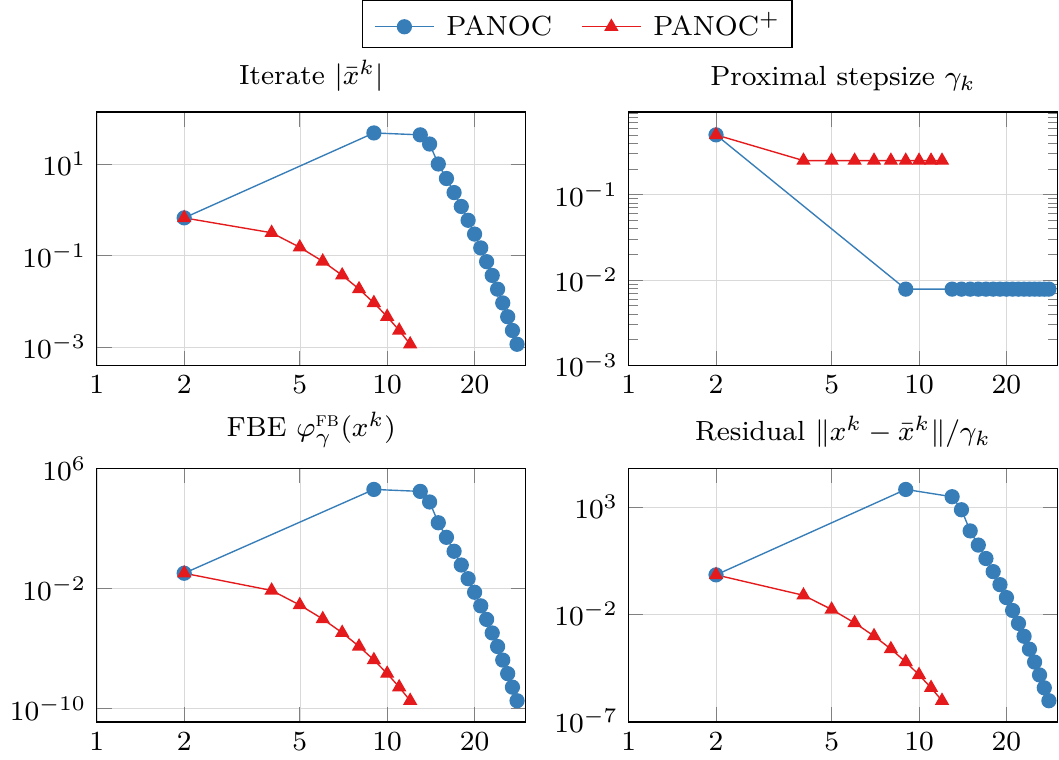}%
		\caption{%
			Comparison of convergence metrics vs number of evaluations of \(\operatorname{T}_\gamma\) for \panoc{} and \panoc{+} on problem \eqref{eq:PANOC_Ex_P_box_smooth}, with Newton-type directions as in \eqref{eq:PANOC_Ex_P_box_smooth_FBEdirection} and parameters \(x^0 = 1\), \(\gamma_0 = 1\), \(\alpha = 0.95\), and \(\beta = 0.5\).
			Similarly to the situation depicted in \cref{fig:counterexamplebox_merit}, the poor geometry estimation of \panoc{} is responsible for an initial divergent behavior that causes slower asymptotic convergence with a small stepsize.%
		}%
		\label{fig:exampleboxnewtonfbe}%
	\end{figure}

	\Cref{fig:exampleboxnewtonfbe} shows that \panoc{}'s iterates initially diverge, even if the starting point $x^0$ is close to the solution $x^\star$, if the proximal stepsize $\gamma_0$ is poorly estimated, in line with the observations above, and despite the choice of regularized Newton-type directions.
	Conversely, \panoc{+} adaptively constructs a tighter merit function and exhibits monotone decrease of \(\FBE\), as depicted in \cref{fig:exampleboxnewtonfbe}.
	Once again, these simulations show that \panoc+ provides a more consistent performance without necessarily requiring more iterations or function evaluations; moreover, the nested linesearch procedure does not lead to a smaller stepsize nor does it hinder fast asymptotic convergence.

%% file: TeX/Text/Analysis.tex
{\def\iter{}%
	In this section we analyze the properties of the iterates generated by \panoc+, starting from their well definedness.
	As a substantial proof of robustness with respect to inexact prox evaluations, we will generalize the setting to an extent that the oracle of the proximal mapping is not required, and instead only a local solution of the proximal sub\-mi\-ni\-mi\-zation problem is needed.
	We will refer to this variant as the \emph{inexact} \panoc{+}, and emphasize that the exact counterpart described in \cref{alg:PANOC+} falls as a special case.

	The investigation in this section originates essentially from three observations.
	Firstly, in the inexact scenario we cannot avail ourselves of the FBE, as its evaluation requires global optimality in the solution of the proximal subproblem.
	Secondly, by considering the equivalent reformulation of \eqref{eq:P}
	\[
		\minimize_{x,z\in\R^n}f(x)+g(z)\quad\stt x=z
	\]
	and defining the associated augmented Lagrangian function
	\begin{align}\label{eq:L}
		\LL_\beta(x,z,y)
	{}\coloneqq{} &
		f(x)+g(z)+\innprod{y}{x-z}+\tfrac\beta2\|x-z\|^2,
	\shortintertext{we remark that}\label{eq:L:FBE}
		\FBE(x)
	{}={} &
		\LL(x,\bar x,-\nabla f(x)),
	\shortintertext{where}\label{eq:L:barx}
		\bar x
	{}\in{} &
		\T(x)
		{}={}
		\argmin\LL(x,{}\cdot{},-\nabla f(x))
	\end{align}
	is the result of an exact proximal minimization.
}%
Thirdly, in the ALM framework, algorithms can be constructed that converge in some sense to stationary points of the optimization problem, even solving the associated subproblems only approximately \cite{birgin2014alm}.
Therefore, we seek relaxed (sub)op\-ti\-mal\-i\-ty concepts for the evaluation of the proximal mapping.
This viewpoint will ultimately highlight how additionally to being used as a solver within ALMs, as in \cite{sopasakis2020open,pas2021alpaqa,demarchi2022constrained}, \panoc{+} can operate as an ALM-type solver itself.

In the broadest possible setting, we do not require any (sub)optimality in the proximal minimization subproblem other than improvement with respect to the previous iteration.
Clearly, additional conditions are needed for generating meaningful iterates, but as a proof of robustness of \panoc+ we demonstrate that any choice complying with said requirement maintains the well definedness of the algorithm.
We will then provide instances of such conditions that, possibly under additional assumptions on the problem, ensure optimality conditions for the limit points of the proposed inexact variant.

Specifically, we consider \cref{alg:PANOC+} with the following instruction replacing \cref{state:PANOC+:barx} therein,
remarking that ``exact'' \(\bar x^k\in\T(x^k)\) as prescribed in \cref{alg:PANOC+} comply with this relaxed requirement (any such \(\bar x^k\) is a global minimizer of \(\LL(x^k,{}\cdot{},-\nabla f(x^k))\), and \(\Phi_k=\FBE(x^k)\) in this case).

\begin{mybox}
	{\it Suboptimal prox step for inexact \panoc+}
	\begin{algorithmic}[1]%
	\makeatletter
		\renewcommand{\alglinenumber}[1]{\footnotesize\fillwidthof[l]{\oldstylenums{88}:}{{\bf 2}.4$'$:}~}%
		\renewcommand{\theALG@line}{2.\oldstylenums{4}$'$}%
	\makeatother
	\State\label{state:NoProx:barx}%
		Let \(\bar x^k\) be a suboptimal minimizer of \(\LL(x^k,{}\cdot{},-\nabla f(x^k))\) such that
		\begin{equation}\label{eq:prox:better}
			\Phi_k
		{}\coloneqq{}
			\LL(x^k,\bar x^k,-\nabla f(x^k))
		{}\leq{}
			\LL(x^k,\bar x^{k-1},-\nabla f(x^k)).
		\end{equation}
	\end{algorithmic}
\end{mybox}

%% file: TeX/Text/Analysis/Results.tex
A crucial complication that the stepsize adjustment in the ``good'' \panoc+ suffers if compared with the original one in the ``bad'' \panoc{}, is that it gives rise to a nested dependency between \(\gamma_k\), \(\tau_k\), and \(d^k\) that could potentially give rise to infinite recursions.
While this is fortunately not the case, as we are about to show, the proof is not as straightforward as in \cite{stella2017simple}.
On top of this, while in the ``exact'' case local boundedness properties of the PG operator \(\T\) could conveniently be exploited, in accounting also for inexactness even for a fixed \(x^k\) the set of points \(\bar x^k\) complying with the relaxed requirement \eqref{eq:prox:better} may be unbounded.
The following result will serve as surrogate of local boundedness for the suboptimal proximal operator.

\begin{lem}\label{thm:inexact:lb}%
	Let a constant \(c\in\R\), a sequence \(\seq{\gamma_j}[j\in\N]\searrow0\), and two bounded sequences \(\seq{u^j,z^j}[j\in\N]\) in \(\R^n\) be fixed, and for every \(j\in\N\) let \(\bar z^j\) be such that
	\[
		g(\bar z^j)+\innprod*{u^j}{\bar z^j-z^j}+\tfrac{1}{2\gamma_j}\|\bar z^j-z^j\|^2
		{}\leq{}
		\tfrac{c}{2\gamma_j}.
	\]
	Then, \(\seq{\bar z^j}[j\in\N]\) is bounded.
	\begin{proof}
		An application of Young's inequality on the inner product yields
		\[
		2\gamma_jg(\bar z^j)
		{}\leq{}
		c
		{}+{}
		\gamma_j\|u_j\|^2
		{}-{}
		(1-\gamma_j)\|\bar z^j-z^j\|^2.
		\]
		To arrive to a contradiction, up to extracting if necessary, suppose that \(0<\|\bar z^j\|\to\infty\).
		Since \(\liminf_{j\to\infty}{g(\bar z^j)}/{\|\bar z^j\|^2}>-\infty\) by \cite[Ex. 1.24]{rockafellar1998variational}, dividing by \(\|\bar z^j\|^2\) and passing to the limit leads to the contradiction \(0\leq -1\).
	\end{proof}
\end{lem}

To avoid trivialities, in what follows we assume that \(x^k\neq\bar x^k\) always holds.
This is consistent with stopping criteria based on the PG residual \(\frac{1}{\gamma_k}\|x^k-\bar x^k\|\), see \cref{sec:Termination}, in which case \(x^k=\bar x^k\) would trigger a successful termination.

\begin{lem}[Well definedness of the ``good'' (inexact) \panoc+]\label{thm:PANOC+:finite}%
	Consider the iterates generated by \cref{alg:PANOC+} with inexact proximal evaluation at \cref{state:PANOC+:barx} as given in \eqref{eq:prox:better}.
	The following hold:
	\begin{enumerate}
	\item\label{thm:finite:LS}%
		Well definedness: at every iteration, the number of backtrackings at \cref{state:PANOC+:gammaLS,state:PANOC+:tauLS} is finite.
	\item\label{thm:finite:descent}%
		At the end of the \(k\)-th iteration (\(k\geq1\)), one has
		\begin{equation}\label{eq:SD}
				\varphi(\bar x^k)
				{}+{}
				\delta_k
			{}\leq{}
				\Phi_k
			{}\leq{}
				\def\iter{k-1}
				\Phi_{k-1}
				{}-{}
				\beta\delta_{k-1}
			\quad\text{where}\quad
				\delta_k
			{}\coloneqq{}
				\tfrac{1-\alpha}{2\gamma_k}
				\|\bar x^k-x^k\|^2.
			\end{equation}
	\item\label{thm:finite:sublevel}%
		Every iterate \(\bar x^k\) remains within \(\lev_{\leq c}\varphi\), where \(c=\Phi_0<\infty\).
	\end{enumerate}
	\begin{proof}
		As observed in \cref{rem:notation}, each iteration \(k\) defines or updates only variables indexed with a \(k\) sub/superscript, while those defined in previous iterations are untouched.
		In what follows, let us index by \({k,j}\) the variables defined at the \(j\)-th attempt within iteration \(k\).
		Note further that \(\gamma_{k,j}L_{k,j}=\alpha\in(0,1)\) holds for every attempt \(j\) within every iteration \(k\), since every time \(\gamma_k\) is halvened the estimate \(L_k\) is doubled (cf. \cref{state:PANOC+:gammaLS}).
		\begin{proofitemize}
		\item\ref{thm:finite:LS}~
			We proceed by induction on \(k\).
			If \(k=0\), there is no backtracking on \(\tau\), and from \cref{thm:inexact:lb} we conclude that all the trials \(\bar x^{0,j}\) remain confined in a bounded set \(\Omega_0\), and therefore any stepsize \(\gamma_{0,j}<\nicefrac{1}{L_{f,\Omega_0}}\) is accepted.

			Suppose now that \(k>0\) and observe that, by the definition of \(\Phi_k\) in \eqref{eq:prox:better} and the failure of the condition at \cref{state:PANOC+:gammaLS}, the inequality
			\begin{equation}\label{eq:lphi}
				\varphi(\bar x^{k-1})
				{}\leq{}
				\Phi_{k-1}
				{}-{}
				\tfrac{1-\alpha}{2\gamma_{k-1}}\|x^{k-1}-\bar x^{k-1}\|^2
			\end{equation}
			holds.
			Since \(\|d^{k,j}\|\leq D\|\bar x^{k-1}-x^{k-1}\|\) and \(\tau_{k,j}\in[0,1]\), any attempt \(x^{k,j}\) defined at \cref{state:PANOC+:x+} during the \(k\)-th iteration satisfies
			\[
				\|x^{k,j}-\bar x^{k-1}\|
			{}={}
				\tau_{k,j}\|
					x^{k-1}-\bar x^{k-1}
					{}+{}
					d^{k,j}
				\|
			{}\leq{}
				(1+D)\|\bar x^{k-1}-x^{k-1}\|
			\]
			and thus remains in a bounded set, be it \(\Omega_k\).
			To arrive to a contradiction, suppose that \(\gamma_{k,j}\searrow0\) as \(j\to\infty\).
			Observe that condition \eqref{eq:prox:better} reads
			\begin{Align*}
				g(\bar x^{k,j})
				{}+{}
				\innprod*{\nabla f(x^{k,j})}{\bar x^{k,j}-\bar x^{k-1}}\quad
			&
			\\
				{}+{}
				\tfrac{1}{2\gamma_{k,j}}\|x^{k,j}-\bar x^{k,j}\|^2
			& {}\leq{}
				g(\bar x^{k-1})
				{}+{}
				\tfrac{1}{2\gamma_{k,j}}\|x^{k,j}-\bar x^{k-1}\|^2.
			\end{Align*}
			Since \(\seq{x^{k,j}}[j\in\N]\) is bounded, an application of \cref{thm:inexact:lb} reveals that \(\seq{\bar x^{k,j}}\) too is bounded.
			Up to possibly enlarging the set, both sequences remain confined in the bounded set \(\Omega_k\), implying that the condition at \cref{state:PANOC+:gammaLS} should have terminated in finite time, whence the sought contradiction.
			
			Hence, \(\gamma_{k,j}\) is backtracked finitely many times within iteration \(k\);
			up to discarding early attempts, we may denote \(\gamma_{k,j}=\gamma_k\).
			Condition \eqref{eq:prox:better} reads
			\begin{align*}
				\LL(x^{k,j},\bar x^{k,j},-\nabla f(x^{k,j}))
			{}\leq{} &
				\LL(x^{k,j},\bar x^{k-1},-\nabla f(x^{k,j}))
			\\
			{}={} &
				f(x^{k,j})+g(\bar x^{k-1})+\innprod*{\nabla f(x^{k,j})}{\bar x^{k-1}-x^{k,j}}
			\\
			&
				{+{}}\tfrac{1}{2\gamma}\|x^{k,j}-\bar x^{k-1}\|^2.
			\end{align*}
			As \(\tau_{k,j}\searrow0\), one has that \(x^{k,j}\to\bar x^{k-1}\).
			Since \(f\) and \(\nabla f\) are continuous, the right-hand side of the inequality converges to \(\varphi(\bar x^{k-1})\), overall resulting in
			\[
			\limsup_{j\to\infty}\LL(x^{k,j},\bar x^{k,j},-\nabla f(x^{k,j}))
			{}\leq{}
			\varphi(\bar x^{k-1})
			{}\overrel[\leq]{\eqref{eq:lphi}}{}
			\Phi_{k-1}-\tfrac{1-\alpha}{2\gamma_{k-1}}\|x^{k-1}-\bar x^{k-1}\|^2.
			\]
			Since \(\|x^{k-1}-\bar x^{k-1}\|>0\) and \(\beta<1\), for \(j\) large enough the condition at \cref{state:PANOC+:tauLS} will be violated and therefore the \(k\)-th iteration successfully terminated.
		\item\ref{thm:finite:descent}
			Follows by combining \eqref{eq:lphi} with the failure of the condition at \cref{state:PANOC+:tauLS} at the end of the iteration.
		\item\ref{thm:finite:sublevel}~
			Direct consequence of assertion \ref{thm:finite:descent}.
		\qedhere
		\end{proofitemize}
	\end{proof}
\end{lem}

We next consider an asymptotic analysis of the algorithm.

\begin{thm}[Asymptotic analysis of the ``good'' (inexact) \panoc+]\label{thm:PANOC+:asymp}%
	Consider the iterates generated by \cref{alg:PANOC+} with inexact proximal evaluation at \cref{state:PANOC+:barx} as given in \eqref{eq:prox:better}.
	The following hold:
	\begin{enumerate}
	\item\label{thm:PANOC+:asymp:cost}%
		\(\seq{\Phi_k}\) converges to a finite value \(\varphi_\star\geq\inf\varphi\) from above.
	\item\label{thm:PANOC+:asymp:summable}%
		\(
			\sum_{k\in\N}\frac{1}{\gamma_k}\|\bar x^k-x^k\|^2
		{}<{}
			\infty
		\).
	\item\label{thm:PANOC+:asymp:fpr}%
		\(
			\lim_{k\to\infty}\|x^k-\bar x^k\|
		{}={}
			\lim_{k\to\infty}\|x^k-x^{k-1}\|
		{}={}
			\lim_{k\to\infty}\|\bar x^k-\bar x^{k-1}\|
		{}={}
			0
		\),
		and in particular the set of limit points of \(\seq{x^k}\) is closed and connected, and coincides with that of \(\seq{\bar x^k}\).
	\item\label{thm:PANOC+:asymp:gamma}%
		\(\sum_{k\in\N}\gamma_k=\infty\).
	\item\label{thm:PANOC+:asymp:res}%
		\(
			\liminf_{k\to\infty}\frac{1}{\gamma_k}\|x^k-\bar x^k\|
		{}={}
			0
		\).
	\item\label{thm:PANOC+:asymp:gammaconstant}%
		Consider the following assertions:
		\begin{enumerate}[label={(\oldstylenums{\arabic*})},ref={(\oldstylenums{\arabic*})}]
		\item\label{thm:PANOC+:gammaconstant:lb}%
			\(\varphi\) is level bounded;
		\item\label{thm:PANOC+:gammaconstant:xk}%
			\(\seq{\bar x^k}\) is bounded;
		\item\label{thm:PANOC+:gammaconstant:barxk}%
			\(\seq{x^k}\) is bounded;
		\item\label{thm:PANOC+:gammaconstant:gamma}%
			\(\seq{\gamma_k}\) is asymptotically constant, \ie, there exists \(\kappa\in\N\) such that \(\gamma_k=\gamma_\kappa\) for every \(k\geq\kappa\);
		\item\label{thm:PANOC+:gammaconstant:Lf}%
			\(f\) has globally Lipschitz-continuous gradient.
		\end{enumerate}
		One has~
		\ref{thm:PANOC+:gammaconstant:lb}
		~\(\Rightarrow\)~
		\ref{thm:PANOC+:gammaconstant:xk}
		~\(\Leftrightarrow\)~
		\ref{thm:PANOC+:gammaconstant:barxk}
		~\(\Rightarrow\)~
		\ref{thm:PANOC+:gammaconstant:gamma}
		~\(\Leftarrow\)~
		\ref{thm:PANOC+:gammaconstant:Lf}.
	\end{enumerate}
	\begin{proof}
		\begin{proofitemize}
		\item\ref{thm:PANOC+:asymp:cost}~
			Follows from \eqref{eq:SD}.
		\item\ref{thm:PANOC+:asymp:summable}~
			A telescoping argument on \eqref{eq:SD} yields
			\begin{equation}\label{eq:summable}
				\beta(1-\alpha)
				\sum_{k\in\N}{
					\tfrac{1}{2\gamma_k}
					\|\bar x^k-x^k\|^2
				}
			{}\leq{}
				\Phi_0
				{}-{}
				\inf\varphi
			{}={}
				\def\iter{0}
				\FBE(x^0)-\inf\varphi,
			\end{equation}
			whence the claimed finite sum.
		\item\ref{thm:PANOC+:asymp:fpr}~
			That \(\|x^k-\bar x^k\|\to0\) follows from assertion \ref{thm:PANOC+:asymp:summable}, since \(\gamma_k\) is upper bounded.
			Next, by the conditions at \cref{state:PANOC+:x+,,state:PANOC+:d}, observe that
			\begin{align*}
			\numberthis\label{eq:xkdiff}
				\|x^k-x^{k-1}\|
			{}={} &
				\bigl\|
					(1-\tau_k)(\bar x^{k-1}-x^{k-1})
					{}+{}
					\tau_kd^k
				\bigr\|
			{}\leq{}
				(1+D)\|\bar x^{k-1}-x^{k-1}\|
			\shortintertext{and thus \(\|x^k-x^{k-1}\|\) vanishes, and in turn so does \(\|\bar x^k-\bar x^{k-1}\|\) since}
				\|\bar x^k-\bar x^{k-1}\|
			{}\leq{} &
				\|x^k-\bar x^k\|
				{}+{}
				\|\bar x^{k-1}-x^{k-1}\|
				{}+{}
				\|x^k-x^{k-1}\|.
			\end{align*}
		\item\ref{thm:PANOC+:asymp:gammaconstant}~
			The first implication follows from \cref{thm:finite:sublevel}, and the second one from assertion \ref{thm:PANOC+:asymp:summable}.
			If \(\seq{x^k}\) is bounded, and thus so is \(\seq{\bar x^k}\), the set \(\Omega_k\) in the proof of \cref{thm:finite:LS} can be taken independent of \(k\), and asymptotic constancy of \(\gamma_k\) follows from the same arguments therein.
			Finally, if \(\nabla f\) is $L_f$-Lipschitz continuous the condition at \cref{state:PANOC+:gammaLS} fails to hold as soon as $\gamma_k\leq\nicefrac{\alpha}{L_f}$ \cite[Prop. A.24]{bertsekas1999nonlinear}, and \(\gamma_k\) is thus asymptotically constant.
		\item\ref{thm:PANOC+:asymp:gamma}~
			By iteratively applying inequality \eqref{eq:xkdiff}, we obtain that
			\begin{align*}
				\|x^k-x^0\|
			{}\leq{} &
				(1+D)\sum_{j=0}^{k-1}\|\bar x^j-x^j\|
			\\
			{}={} &
				(1+D)\sum_{j=0}^{k-1}\gamma_j^{-\nicefrac12}\|\bar x^j-x^j\|\gamma_j^{\nicefrac12}
			\\
			{}\leq{} &
				\textstyle
				(1+D)
				\sqrt{
					\sum_{j=0}^{k-1}
					\gamma_j^{-1}\|\bar x^j-x^j\|^2
				}
				\sqrt{
					\sum_{j=0}^{k-1}\gamma_j
				}
			\\
			{}\overrel*[\leq]{\eqref{eq:summable}}{} &
				\textstyle
				(1+D)
				\sqrt{
					2
					\frac{\def\iter{0}\FBE(x^0)-\inf\varphi}{\beta(1-\alpha)}
				}
				\sqrt{
					\sum_{j=0}^{k-1}\gamma_j
				}.
			\end{align*}
			Contrary to the claim, if \(\sum_{k\in\N}\gamma_k<\infty\) holds, then \(\seq{x^k}\) is bounded.
			From assertion \ref{thm:PANOC+:asymp:gammaconstant} proven above we then infer that \(\gamma_k\) is asymptotically constant, thus contradicting the finiteness of \(\sum_{k\in\N}\gamma_k\).
		\item\ref{thm:PANOC+:asymp:res}~
			Immediate consequence of assertions \ref{thm:PANOC+:asymp:summable} and \ref{thm:PANOC+:asymp:gamma}.
		\qedhere
		\end{proofitemize}
	\end{proof}
\end{thm}

\begin{rem}
	If the iterates remain bounded (as is the case when the objective \(\varphi\) is level bounded), owing to \cref{thm:PANOC+:asymp:gammaconstant}, \cref{alg:PANOC+} with exact prox evaluations as in \cref{state:PANOC+:barx} eventually reduces to the original \panoc{} \cite{stella2017simple} with constant stepsize, and its convergence results are then readily available, including global convergence (possibly at R-linear rates) under Kurdyka-{\L}ojasiewicz assumptions, and superlinear when converging to a strong local minimum with directions satisfying the Dennis-Mor{\'e} condition, see \cite{stella2017simple,themelis2018forward}.
\end{rem}

Nevertheless, even in accounting for inexact proximal evaluations it is still possible to derive some qualitative guarantees for the limit points, provided that \(\bar x^k\) satisfies some local suboptimality requirements.
We list two such instances in the following definition and later detail a proof validating the claim.

\begin{defin}[Prox suboptimality criteria]%
	Relative to the minimization problem \eqref{eq:L:barx} defining the PG mapping, we say that the iterates \(\bar x^k\) computed at \cref{state:NoProx:barx} are:
	\begin{enumerate}
	\item\label{defin:delta}%
		\DEF{\(\delta\)-stationary} (for some \(\delta>0\)) if
		\(\mathtight[0.4]
			\dist\bigl(0,
				\partial\bigl[\LL(x^k,{}\cdot{},-\nabla f(x^k))\bigr](\bar x^k)
			\bigr)
		{}\leq{}
			\delta
		\),
		that is, if there exists \(\bar v^k\in\partial g(\bar x^k)\) such that
		\begin{equation}\label{eq:prox:stationary}
			\bigl\|
				\bar v^k+\nabla f(x^k)+\tfrac{1}{\gamma_k}(\bar x^k-x^k)
			\bigr\|
		{}\leq{}
			\delta.
		\end{equation}
	\item\label{defin:uniform}%
		\DEF{Uniformly locally optimal} if there exist \(r>0\) and a sequence \(\varepsilon_k\searrow0\) such that the following local minimality condition holds:
		\begin{equation}\label{eq:locmin}
			\LL(x^k,\bar x^k,-\nabla f(x^k))
		{}\leq{}
			\LL(x^k,x,-\nabla f(x^k))+\varepsilon_k
		\quad
			\forall x\in\cball{\bar x^k}{r}.
		\end{equation}
	\end{enumerate}
\end{defin}

Notice that no (approximate) local minimality is required in the approximate stationarity criterion of \cref{defin:delta}.
Consequently, the output can be retrieved by any descent method starting at the previous iteration and terminating when \(\delta\)-stationarity is achieved.
It is also worth remarking that the prox suboptimality tolerance \(\delta\) does not need to be small nor fixed for all iterations, and can instead be replaced by a sequence \(\delta_k\searrow\delta\geq0\).
The uniform local optimality requirement of \cref{defin:uniform} is instead more restrictive, and is possibly subject to prior knowledge on the geometry of the augmented Lagrangian.
The uniformity is dictated by the value of \(r>0\), whose role can be appreciated by considering the sequence \(z^k=\nicefrac1k\) for \(k>0\) which consists of (isolated) local minimizers for the function
\[
	h(x)
{}={}
	\begin{ifcases}
		x & x=\nicefrac1k,~k\in\N_{>0}\\
		x^2+x-1 & x\leq 0\\
		\infty\otherwise,
	\end{ifcases}
\]
yet the limit \(z=0\) is not stationary for \(h\).
The pathology arises from the non uniformity of the radius of local minimality of \(z^k\), which is \(r_k<\nicefrac{1}{k(k+1)}\to0\).

\begin{thm}[Subsequential convergence of inexact \panoc+]\label{thm:PANOC+:subseq}%
	Consider the iterates generated by \cref{alg:PANOC+} with inexact proximal evaluation at \cref{state:PANOC+:barx} as given in \eqref{eq:prox:better}.
	Suppose that the iterates remain bounded (as is the case when \(\varphi\) is coercive), and let \(\omega\) be the set of limit points of \(\seq{\bar x^k}\).
	Then:
	\begin{enumerate}
	\item\label{thm:NoProx:subseq:omega:subdiffC0}%
		If \(\seq{\bar x^k}\) are \(\delta\)-stationary as in \cref{defin:delta} and \(\graph\partial g\) is closed relative to \(\dom g\times\R^n\), then \(\omega\) is made of \(\delta\)-stationary points for \(\varphi\).
	\item\label{thm:NoProx:subseq:omega:locmin}%
		If the sequence \(\seq{\bar x^k}\) is (eventually) uniformly locally optimal as in \cref{defin:uniform} (this being true in case of exact prox evaluations, having \(r=\infty\) and \(\varepsilon_k=0\) in this case), then the set \(\omega\) is made of stationary points for \(\varphi\), and \(\varphi\) is constantly equal to \(\varphi_\star\) as in assertion \ref{thm:PANOC+:asymp:cost} there.
	\end{enumerate}
	\begin{proof}
		\def\iter{}%
		Up to possibly discarding early iterates, in light of the boundedness of the sequences and the consequent eventual constancy of \(\gamma_k\) by \cref{thm:PANOC+:asymp:gammaconstant}, we may assume that \(\gamma_k\equiv\gamma>0\) holds for all \(k\).
		Let \(x^\star\in\omega\) be fixed, and let an infinite set of indices \(K\subseteq\N\) be such that \(\seq{\bar x^k}[k\in K]\to x^\star\), so that \(\seq{x^k}[k\in K]\to x^\star\) too as it follows from \cref{thm:PANOC+:asymp:fpr}.
		\begin{proofitemize}
		\item\ref{thm:NoProx:subseq:omega:subdiffC0}~
			Since \(\nabla f(x^k)+\tfrac1\gamma(\bar x^k-x^k)\to\nabla f(x^\star)\) as \(K\ni k\to\infty\), up to extracting a subsequence if necessary, it follows from \eqref{eq:prox:stationary} that \(\bar v^k\to\bar v^\star\) with \(\|\bar v^\star+\nabla f(x^\star)\|\leq\delta\).
			Since \(\seq{\Phi_k=\LL(x^k,\bar x^k,-\nabla f(x^k))}\) is bounded, owing to \cref{thm:PANOC+:asymp:cost}, and since both \(f\) and \(\nabla f\) are continuous, clearly \(\seq{g(\bar x^k)}\) remains bounded, and therefore, by lower semicontinuity, \(x^\star\in\dom g\).
			Since also \(\seq{\bar x^k}[k\in K]\subseteq\dom g\), from the assumptions we conclude that \(\bar v^\star\in\partial g(x^\star)\) and thus \(\bar v^\star+\nabla f(x^\star)\in\partial\varphi(x^\star)\), proving \(\delta\)-stationarity of \(x^\star\) for \(\varphi\).
		\item\ref{thm:NoProx:subseq:omega:locmin}~
			Letting \(\varphi_\star\) be as in \cref{thm:PANOC+:asymp:cost} and invoking \eqref{eq:SD}, lsc of \(\varphi\) yields \(\varphi(x^\star)\leq\varphi_\star\).
			For \(k\) large enough so that \(\bar x^k\) is \(r\)-close to \(x^\star\), we have
			\begin{align*}
				\varphi_\star
			{}={}
				\lim_{k\in K}
				\Phi_k
			{}={} &
				\lim_{k\in K}
				\LL(x^k,\bar x^k,-\nabla f(x^k))
			\\
			{}\leq{} &
				\limsup_{k\in K}
				\LL(x^k,x^\star,-\nabla f(x^k))+\varepsilon_k
			\\
			{}={} &
				\LL(x^\star,x^\star,-\nabla f(x^\star))
			{}={}
				\varphi(x^\star)
			{}\leq{}
				\varphi_\star,
			\end{align*}
			owing to continuity of \(f\) and \(\nabla f\), and the fact that both \(\varepsilon_k\) and \(\|x^k-\bar x^k\|\) vanish (the former by assumption and the latter by \cref{thm:PANOC+:asymp:fpr}).
			From the arbitrarity of \(x^\star\in\omega\) we conclude that \(\varphi\) is constant on \(\omega\) with value \(\varphi_\star\).
			Notice further this also shows that \(g(\bar x^k)\to g(x^\star)\) as \(K\ni k\to\infty\).
			Ekeland's variational principle \cite[Prop. 1.43]{rockafellar1998variational} with \(\delta_k=\sqrt{\varepsilon_k}\) ensures for every \(k\in K\) (large enough so that \(\sqrt{\varepsilon_k}\leq r\)) the existence of \(\xi^k\in\cball{\bar x^k}{\sqrt{\varepsilon_k}}\)
			together with
			\[
				\eta^k
			{}\in{}
				\hat\partial\bigl[\LL(x^k,{}\cdot{},-\nabla f(x^k))\bigr](\xi^k)
			{}={}
				\nabla f(x^k)+\hat\partial g(\xi^k)+\tfrac{1}{\gamma}(\xi^k-x^k)
			\]
			such that
			\(
				\LL(x^k,\xi^k,-\nabla f(x^k))
			{}\leq{}
				\Phi_k
			\)
			and \(\eta^k\in\cball{0}{\sqrt{\varepsilon_k}}\).
			By lsc of \(g\) and since \(\xi^k\to x^\star\), necessarily \(g(\xi^k)\to g(x^\star)\) and the inclusion \(-\nabla f(x^\star)\in\partial g(x^\star)\) is then readily obtained, whence the claimed stationarity of \(x^\star\) for \(\varphi\).
		\qedhere
		\end{proofitemize}
	\end{proof}
\end{thm}

Closedness of \(\graph\partial g\) relative to \(\dom g\times\R^n\) as required in \cref{thm:NoProx:subseq:omega:subdiffC0} is frequently encountered in applications, and trivially encompasses all functions that are continuous on their domain, such as indicators of closed sets.
The 0-norm is instead an example of a function which is not continuous on its domain but that nevertheless complies with the requirement in \cref{thm:NoProx:subseq:omega:subdiffC0}.
Indeed, notice that
\[
	\partial g(x)
{}={}
	\hat\partial g(x)
{}={}
	E_1\times\dots\times E_n,
\quad\text{where}\quad
	E_i
{}={}
	\begin{ifcases}
		\R & x_i=0\\
		\set0 & x_i\neq 0
	\end{ifcases}
\]
for \(g=\|{}\cdot{}\|_0\).
Consider a sequence \(x^k\to x\) along with \(\partial g(x^k)\ni v^k\to v\); we will show that \(v\in\partial g(x)\), regardless of whether or not \(g(x^k)\) converges to \(g(x)\).
Indeed, if \(x_i=0\), then trivially \(v_i\in\R=E_i\).
Otherwise, \(x_i^k\neq 0\) holds for large enough \(k\), thus necessarily \(v_i^k=0\), and consequently \(v_i\in\set{0}=E_i\).
Either way, since this holds for every component, we conclude that \(v\in\partial g(x)\).

%% file: TeX/Text/Analysis/Termination.tex
\Cref{alg:PANOC+} runs indefinitely and generates an infinite sequence of iterates \(\seq{x^k}\) and \(\seq{\bar x^k}\).
Along its execution, we are compelled to check some suitable conditions for stopping and returning an \(\bar x^k\) that, in some sense, satisfactorily minimizes \(\varphi\).
The assertion of \cref{thm:PANOC+:asymp:res} guarantees that the standard termination criterion on the residual
\begin{equation}\label{eq:termination}
	\tfrac{1}{\gamma_k}\|x^k-\bar x^k\|\leq\tfrac{\varepsilon}{2}
\end{equation}
is verified in finite time.
However, considering \eqref{eq:proxoptim}, a control on the magnitude of \(\|\nabla f(x^k)-\nabla f(\bar x^k)\|\) must also be imposed in order to guarantee bounds on \(\dist(0,\partial\varphi(\bar x^k))\).
This calls for a strengthened linesearch condition at \cref{state:PANOC+:gammaLS} ensuring also the satisfaction of
\begin{equation}\label{eq:nabladiff}
	\|\nabla f(x^k)-\nabla f(\bar x^k)\|
{}\leq{}
	\tfrac{1}{\gamma_k}\|x^k-\bar x^k\|,
\end{equation}
so that, by a triangular inequality argument on \eqref{eq:proxoptim}, \(\varepsilon\)-stationarity of \(\bar x^k\) (that is, \(\dist(0,\partial\varphi(\bar x^k))\leq\varepsilon\)) would be guaranteed by \eqref{eq:termination}.
On the one hand, owing to \cref{ass:f} the proof of \cref{thm:finite:LS} (and of all other results) would still verbatim apply, meaning that this criterion would not affect the well definedness of \cref{alg:PANOC+}, or in fact any result presented so far.
On the other hand, this would require evaluations of \(\nabla f(\bar x^k)\), otherwise not needed, and thus affect the overall complexity.
To account for this fact, a viable solution is to trigger this strengthened linesearch only after \eqref{eq:termination} is first satisfied, at which point the algorithm can terminate whenever \eqref{eq:termination} is verified again.

Note that the same conclusions can be made under suboptimal prox evaluations complying with the local uniformly of \cref{defin:uniform}, as long as \(\varepsilon_k=0\) for all \(k\).
In case of \(\delta\)-stationarity as in \cref{defin:delta}, instead, the same criterion would guarantee \((\delta+\varepsilon)\)-stationarity of the output.

%% file: TeX/Text/Analysis/NM.tex
Nonmonotone linesearch procedures often prove beneficial in practice, as they can reduce conservatism in the linesearch and favor larger steps.
By patterning the rationale of the ZeroFPR algorithm \cite{themelis2018forward}, a nonmonotone linesearch can be readily integrated in \panoc+ at \cref{state:PANOC+:tauLS} without affecting the finite termination and asymptotic properties asserted in \cref{thm:PANOC+:finite,thm:PANOC+:asymp}.
This is done by changing the definition of \(\Phi_k\) at \cref{state:PANOC+:barx} into \(\Phi_k=(1-p_k)\Phi_{k-1}+p_k\FBE(x^k)\) for \(k>0\) (with \(\FBE(x^k)\) being replaced by \(\LL(x^k,\bar x^k,-\nabla f(x^k))\) in the inexact case), where \(\seq{p_k}\subset(0,1]\) is any user-selected sequence bounded away from 0.
The key observation enabling the possibility to replicate all the convergence results is the inequality \(\FBE(x^k)\leq\Phi_k\), which follows from an elementary induction (cf. \cite[Lem. 5.1]{themelis2018forward}).

%% file: TeX/Text/Analysis/PG.tex
By selecting \(d^k = \bar x^{k-1} - x^{k-1}\) at \cref{state:PANOC+:d}, \panoc+ reduces to the classical proximal gradient method \(\def\iter{k-1}x^k\in\T(x^{k-1})\) with an adaptive stepsize.
In fact, the descent condition at \cref{state:PANOC+:tauLS} does not need to be checked, as it is always satisfied for any \(\tau_k\), having
\(
	x^k
{}={}
	(1-\tau_k)\bar x^{k-1}+\tau_k(x^k+d^k)
{}={}
	\bar x^{k-1}
\)
independently of the value of \(\tau_k\).
For this specific choice of the update direction \(d^k\), the algorithm simplifies and reduces to the proximal gradient method with adaptive stepsize selection given in \cref{alg:PG}.
Convergence results developed in the general setting of \panoc+ can thus be readily imported, even in the inexact case.

\begin{cor}[Convergence of adaptive PG]%
	All the assertions of \cref{thm:PANOC+:asymp,thm:PANOC+:subseq} remain valid for the iterates generated by \cref{alg:PG}.
\end{cor}

\begin{algorithm}[h]
	\caption{Inexact proximal gradient with adaptive \(\gamma\)-stepsize rule}%
	\label{alg:PG}%
	\input{TeX/Alg/PG.tex}\end{algorithm}

We note that the exact version of \cref{alg:PG}, that is, with \(\bar x^k\in\T(x^k)\) in \cref{state:PG:barx}, corresponds to a simplified version of the linesearch strategy \cite[LS1]{salzo2017variable}, with no relaxation and in finite dimensional spaces but here analyzed for (fully) nonconvex problems.
Alternatively, it can be viewed as the monotone PG method outlined in \cite[Alg. 3.1]{kanzow2021convergence} with a slightly more conservative linesearch, since
\begin{align*}
		\varphi(\bar x^k)
	{}\leq{} &
		f(x^k)
		{}+{}
		\innprod{\nabla f(x^k)}{\bar x^k-x^k}
		{}+{}
		\tfrac{\alpha}{2\gamma_k}\|\bar x^k-x^k\|^2
		{}+{}
		g(\bar x^k)
	\\
	{}\overrel*[=]{\eqref{eq:FBEmin}}{} &
		\FBE(x^k) - \tfrac{1-\alpha}{2\gamma_k}\|\bar x^k-x^k\|^2
	{}\leq{}
		\varphi(x^k) - \tfrac{1-\alpha}{2\gamma_k}\|\bar x^k-x^k\|^2,
\end{align*}
where the inequalities follow from \cref{state:PG:gammaLS} and \cref{thm:FBE:leq}.
Remarkably, plain continuous differentiability (as opposed to locally Lipschitzian) suffices in the given reference, under a few other technical assumptions.
However, the discussion therein is confined to plain PG iterations as in \cref{alg:PG}, while our analysis is more general and captures plain PG as simple byproduct.

%% file: TeX/Alg/PG.tex
\begin{algorithmic}[1]%
\linespread{1.4}\selectfont%
\Require
	\(x^0\in\R^n\);~~
	\(\gamma_0\in(0,\gamma_g)\);~~
	\(\alpha\in(0,1)\)%
\Initialize
	\(\bar x^{-1}=x^0\),~~
	\(k\gets0\),~~
	and start from \cref{state:PG:barx}
\vspace*{1ex}%
\hrule
\vspace*{1ex}%
\def\iter{k}%
	\State
		\label{state:PG:init}%
		\(\gamma_k\gets\gamma_{k-1}\),~
		\(x^k\gets\bar x^{k-1}\)
	\State
		\label{state:PG:barx}%
		Let \(\bar x^k\) be as in \eqref{eq:prox:better}
		(\eg, \(\bar x^k\in\T(x^k)\))
	\If{
		\(
			f(\bar x^k)
		{}>{}
			f(x^k)
			{}+{}
			\innprod{\nabla f(x^k)}{\bar x^k-x^k}
			{}+{}
			\tfrac{\alpha}{2\gamma_k}\|\bar x^k-x^k\|^2
		\)
	}%
	\label{state:PG:gammaLS}%
		\Statex*%
			\(\gamma_k\gets\nicefrac{\gamma_k}{2}\),~
			and go back to \cref{state:PG:barx}
	\EndIf
	\State
		\(k\gets k+1\)~
		and start the next iteration at \cref{state:PG:init}
\end{algorithmic}

%% file: TeX/Text/Conclusions.tex
We investigated an adaptive scheme to appropriately select the proximal stepsize within solvers for fully nonconvex composite optimization, focusing on (and extending) the PANOC framework.
Our convergence analysis demonstrates the well-definedness of the algorithm and characterizes its asymptotic properties, possibly in the absence of (global) Lipschitz gradient continuity for the smooth term.
Indeed, witnessing the approach's robustness, we considered a setting with possibly inexact proximal mapping oracle for the nonsmooth term, providing suitable conditions for its approximate computation.
By means of detailed illustrative examples, we highlighted weaknesses of previous approaches and the crucial steps undertaken in this work, as well as their benefits in terms of convergence guarantees and efficiency.
Our findings indicate that, by better capturing the problem's geometry, a more conservative adaptive scheme can yield superior practical performance under weaker conditions.
Comprising also arbitrary acceleration directions and nonmonotone variants, these results significantly enlarge the scope of PANOC, both as stand-alone tool for optimization and internal solver within other algorithms, \eg in ALM and sequential programming approaches.